\documentclass{article}
\usepackage{graphicx} 
\usepackage{amsmath,amssymb,amsthm}
\usepackage[margin=1.1in]{geometry}
\usepackage{xcolor}
\usepackage{tikz}
\usetikzlibrary{decorations.markings,arrows.meta,braids}

\usepackage[colorlinks=true,linkcolor=blue!50!black,citecolor=blue!50!black,urlcolor=blue!50!black]{hyperref}

\newtheorem{theorem}{Theorem}
\newtheorem{proposition}[theorem]{Proposition}
\newtheorem{lemma}[theorem]{Lemma}
\newtheorem{corollary}[theorem]{Corollary}
\theoremstyle{definition}
\newtheorem{definition}[theorem]{Definition}
\newtheorem{example}[theorem]{Example}
\newtheorem{remark}[theorem]{Remark}
\newtheorem{observation}[theorem]{Observation}

\title{There are hard braids}
\author{Alexei Vernitski}
\date{November 2026}

\definecolor{topc}{RGB}{190,60,40}
\definecolor{botc}{RGB}{40,90,170}
\newcommand{\fl}[1]{{\setlength{\fboxsep}{1.2pt}\colorbox{orange!22}{$#1$}}}
\begin{document}

\maketitle

\begin{abstract}
Let us call a trivial braid hard if it cannot be untangled by Reidemeister moves without increasing the number of crossings. Are there hard trivial braids? Our paper answers this question.

Hard braids exist. For every number of strands starting from 4, for every even length starting from 16, we present examples of hard trivial braids of that size.

Our smallest example is minimal. We prove that there are no hard braids on three strands.
On four strands our shortest example, of length 16, is a shortest hard braid.

The hard braids in that family can be untangled without increasing their length if we use Reidemeister moves and an infinite family of moves
$\sigma_i^{\varepsilon}\sigma_j^{m}\sigma_i^{-\varepsilon}
\leftrightarrow\sigma_j^{-\varepsilon}\sigma_i^{m}\sigma_j^{\varepsilon}$, which involve only
three strands.

We also build a harder family of hard braids which no moves confined to a bounded number of strands can untangle. For every $m$ we construct a trivial braid on $4m$ strands, of length $16m^2$, which can be untangled without increasing its length by moves on $2m+1$ strands but not not on $2m$ strands.
\end{abstract}

\section{Introduction}

Can one untangle a trivial braid without increasing the number of crossings? This
is the question we study in this article. The answer depends on the moves one allows, and
we begin with the most basic ones.

Throughout, $B_n$ is the braid group on $n$ strands with Artin's presentation
\[
  B_n=\bigl\langle \sigma_1,\dots,\sigma_{n-1}\ \bigm|\
  \sigma_i\sigma_j=\sigma_j\sigma_i\ (|i-j|\ge 2),\quad
  \sigma_i\sigma_j\sigma_i=\sigma_j\sigma_i\sigma_j\ (|i-j|=1)\bigr\rangle .
\]
A braid word is a word in the letters $\sigma_i^{\pm1}$; its length is its number
of crossings. A braid word is \emph{trivial} if it represents $1$ in $B_n$.

\begin{definition}
The \emph{elementary moves} on a braid word are:
\begin{enumerate}
  \item deleting a subword $\sigma_i^{\varepsilon}\sigma_i^{-\varepsilon}$;
  \item replacing $\sigma_i^{\varepsilon}\sigma_j^{\eta}$ by
    $\sigma_j^{\eta}\sigma_i^{\varepsilon}$ when $|i-j|\ge2$;
  \item replacing $\sigma_i^{x}\sigma_j^{y}\sigma_i^{z}$ by
    $\sigma_j^{z}\sigma_i^{y}\sigma_j^{x}$ when $|i-j|=1$, where $x,y,z\in\{1,-1\}$ and the
    exponents do not alternate, that is, unless $x=z\ne y$.
\end{enumerate}
Moves 2 and 3 are exactly the replacements of one half of a defining relator by
the other half: the patterns $\sigma_i^{\pm1}\sigma_j^{\mp1}\sigma_i^{\pm1}$ never
occur as halves of relators. None of the moves increases the length. A trivial
braid word is \emph{hard} if it cannot be transformed into the empty word by
elementary moves, even when cyclic permutations of the word and the passage from a
word to its reversed inverse are also allowed.
\end{definition}

We first recall the context of this question.

\paragraph{Words in Coxeter groups.} In combinatorial group theory, Tits' theorem
deals with a class of groups closely related to braid groups. A Coxeter group is
presented by involutions $s_i$ subject to braid-type relations; the symmetric group
$\mathfrak S_n$, obtained from $B_n$ by adding the relations $\sigma_i^2=1$, is the
basic example. Tits proved that any word in the generators of a Coxeter group can
be transformed into a reduced word, of minimal length, by braid relations and
deletions of factors $s_is_i$, and that any two reduced words for the same
element are related by braid relations alone \cite{Tits1969}; the second
statement was obtained independently by Matsumoto \cite{Matsumoto1964}. In
particular, a word representing the identity can be transformed into the empty
word by a sequence of steps which never increase its length. Standard modern
accounts are \cite[Ch.~IV]{Bourbaki} and \cite[\S3.3]{BjornerBrenti}.

For the groups that interest us here, Artin--Tits groups such as $B_n$, the closest
analogue is the conjecture of Dehornoy and Godelle that the word problem can be
solved without introducing trivial factors $ss^{-1}$ or $s^{-1}s$, their
Property~H \cite{DehornoyGodelle}. They proved it for Artin--Tits groups of type
FC, a class that includes the braid groups. Property~H forbids the insertion of
trivial factors but, unlike Tits' theorem, places no bound on the length of the
intermediate words; it is therefore weaker than the absence of hard braids.

\paragraph{Diagrams of knots.} In knot theory, it is known that some diagrams of
the trivial knot are \emph{hard}, meaning that when we simplify them using
Reidemeister moves, it is necessary to increase the number of crossings in the
diagram before we can decrease it. The classical example goes back to Goeritz
\cite{Goeritz}; collections of hard unknot diagrams, and measures of their
hardness, are discussed by Kauffman and Lambropoulou \cite{KauffmanLambropoulou}
and by Burton, Chang, L\"offler, de~Mesmay, Maria, Schleimer, Sedgwick and Spreer
\cite{BurtonEtAl}, who exhibit diagrams requiring at least three extra crossings;
hard diagrams of split links are studied in \cite{LunelDeMesmaySpreer}, and
millions of hard unknot diagrams were found with reinforcement learning in
\cite{ApplebaumEtAl}. On the other hand, the number of Reidemeister moves needed is
bounded by a polynomial in the number of crossings \cite{Lackenby}, improving the
exponential bound of Hass and Lagarias \cite{HassLagarias}; quadratically many
moves may be necessary \cite{HassNowik}.

Two results show that hardness depends on the type of diagram and on the moves
allowed. Dynnikov proved that rectangular diagrams (arc-presentations) of the
unknot \emph{can} be simplified monotonically, by exchange moves and
destabilisations that never increase the complexity \cite{Dynnikov2006}. For
closed braids, Morton found a closed $4$-braid with unknotted closure that cannot
be reduced by conjugation and destabilisation without first increasing the number
of strands \cite{Morton1983}; Birman and Menasco proved that closed braid
representatives of the unlink can be simplified to the trivial closed braid using
braid isotopy, destabilisations and exchange moves, without increasing the braid
index \cite{BirmanMenasco}.

It is from this literature that we borrow the word \emph{hard}. Summarising the context
above: if braids behave like words in Coxeter groups, there should be no hard braids; if
they behave like diagrams of knots, there should be some.

\paragraph{Algorithmic complexity.} The word problem in $B_n$ has a long history.
Garside solved it in 1969 \cite{Garside}; the automatic structure of braid groups
\cite{Epstein} and the Garside normal form \cite{ElrifaiMorton} give algorithms
that are quadratic in the length of the word for a fixed number of strands, and the
algorithm of Birman, Ko and Lee \cite{BirmanKoLee} runs in time $O(\ell^2 n)$ for
words of length $\ell$ in $B_n$, a bound also obtained by Hamidi-Tehrani
\cite{HamidiTehrani}. Dehornoy's handle reduction \cite{Dehornoy1997} and Dynnikov's
coordinates \cite{Dynnikov2002} are very fast in practice; see
\cite{DDRW} for both. With Dynnikov's coordinates, comparing two braid words of
length $\ell$ takes time $O(\ell^2)$, since the coordinates have $O(\ell)$ digits.
Very recently Bell and Schleimer announced an $O(\ell\log^3\ell)$ algorithm for the
word problem in mapping class groups of compact surfaces \cite{BellSchleimer},
building on Dynnikov's work on normal curves \cite{Dynnikov2022} and on fast integer
arithmetic \cite{HarveyVanDerHoeven}.

If there were no hard braids, the word problem in $B_n$ would admit a strikingly
simple solution: a word is trivial if and only if the empty word can be reached from
it by elementary moves, so that only words no longer than the input need ever be
considered. Such a procedure need not be fast, since the words reachable by
length-preserving moves can be numerous; but it would be a purely combinatorial
characterisation of trivial braids of the kind Tits' theorem provides for Coxeter
groups. The results below show that no such characterisation holds, neither for the
elementary moves nor for any family of moves confined to a bounded number of strands.

\paragraph{Our contributions.} This article makes four contributions.
\begin{enumerate}
\item \emph{Hard braids exist, and there are many of them.} Theorem~\ref{thm:hard}
  (Section~\ref{sec:family}) gives, for every $n\ge4$, an explicit infinite family of
  hard braid words $w_\varepsilon(p,q)$ in $B_n$, depending on the position $j$ of three
  consecutive generators, a sign $\varepsilon$, and two integers $p,q$ with $|p|,|q|\ge2$.
  Its members have length $2|p|+2|q|+8$, so there are hard braids of every even length
  from $16$ on, in every $B_n$ with $n\ge4$. The shortest member is a word $W\in B_4$ of
  length $16$, which Section~\ref{sec:diagram} draws as a van Kampen diagram.
\item \emph{Our hard braids are minimal in the number of strands.} Theorem~\ref{thm:three}
  (Section~\ref{sec:three}) shows that there are no hard braids on three strands; on two
  strands there are none for trivial reasons. Hence four strands are needed, and our
  family attains this bound. The proof is a complete combinatorial argument with van
  Kampen diagrams, and it identifies the feature of four strands that hard braids exploit.
  An exhaustive computer search (Section~\ref{sec:verification}) shows moreover that on four
  strands no hard braid is shorter than $W$.
\item \emph{The hard braids of this family are not very hard.} Theorem~\ref{thm:power}
  (Section~\ref{sec:power}) shows that every member of the family can be untangled
  without ever increasing its length once power braid moves are allowed; exactly $2|p|$
  of them are needed. Power braid moves are themselves length-preserving consequences of
  the braid relations, so the hardness of our braids does not force any extra crossings:
  it only shows that the elementary moves, which rearrange at most three crossings at a
  time, are too local. How non-local the moves must be is, however, unbounded:
  Proposition~\ref{prop:exponent} shows that the braids $W_k$ obtained from $W$ by
  replacing every $\sigma_2^{\pm2}$ by $\sigma_2^{\pm2k}$ need power braid moves of
  exponent $2k$, so no bound on the exponent of power braid moves suffices for all hard
  braids. Like elementary braid moves, power braid moves involve only three strands at a
  time (Remark~\ref{rem:threestrands}); what the braids of this family require is
  therefore moves that are local in the strands but not in the number of crossings.
\item \emph{No moves confined to a bounded number of strands untangle all trivial braids.}
  The third contribution suggests an analogue of Tits' theorem for braid groups: perhaps
  every trivial braid can be untangled monotonically by moves, of arbitrary length, each of
  which takes place inside a subgroup generated by two standard generators, and so
  involves at most three strands. Section~\ref{sec:bands} shows that this is false, and that
  no bound on the number of strands can replace three. A \emph{move on $k$ strands} replaces
  a part of a braid word in which only $k$ adjacent strands take part by an equal braid on
  the same strands with no more crossings. For every $m\ge1$, Theorem~\ref{thm:bands} gives
  a trivial braid word $V_m$ on $4m$ strands, with $16m^2$ crossings, such that no move on
  at most $2m$ strands does anything to $V_m$ beyond commuting distant crossings, whereas
  moves on $2m+1$ strands untangle it. The word $V_m$ is obtained from a hard braid of our
  family by replacing each strand by a band of $m$ parallel strands; $V_1$ is that hard
  braid, and $V_2$, with $64$ crossings on $8$ strands, is our counterexample to the
  analogue of Tits' theorem. The proof of rigidity is combinatorial and is carried out by
  hand for all $m$; it reduces the claim to a few properties of the hard braid $V_1$. For
  $m=2$ and $m=3$ it is confirmed by exhaustive computer searches.
\end{enumerate}
Section~\ref{sec:verification} describes the computer checks that accompany the proofs,
and Section~\ref{sec:conclusion} summarises the results and lists open questions.

\section{An infinite family of hard braids}\label{sec:family}

This section establishes our first contribution: hard braids exist, and there are many
of them. The main result is Theorem~\ref{thm:hard}.

\subsection{Power braid moves}

Besides the elementary moves we shall use the following length-preserving
replacements.

\begin{definition}
For $|i-k|=1$, $\varepsilon=\pm1$ and $m\in\mathbb Z\setminus\{0\}$, the
\emph{power braid move} is
\[
  \sigma_i^{\varepsilon}\,\sigma_k^{m}\,\sigma_i^{-\varepsilon}
  \ \longleftrightarrow\
  \sigma_k^{-\varepsilon}\,\sigma_i^{m}\,\sigma_k^{\varepsilon}.
\]
For $m=1$ this is the elementary braid move
$\sigma_i^{\varepsilon}\sigma_k\sigma_i^{-\varepsilon}\leftrightarrow
\sigma_k^{-\varepsilon}\sigma_i\sigma_k^{\varepsilon}$, whose two sides are conjugates of
$\sigma_k$ and of $\sigma_i$ representing the same braid. Raising both sides to the $m$-th
power, and cancelling the factors $\sigma_i^{-\varepsilon}\sigma_i^{\varepsilon}$ and
$\sigma_k^{\varepsilon}\sigma_k^{-\varepsilon}$ that appear, shows that the two sides of the
power braid move represent the same braid for every $m$. They have the same length,
$|m|+2$.
\end{definition}

Like an elementary braid move, a power braid move involves only three strands: both
sides are words in two adjacent generators, which act on three consecutive positions.
What distinguishes it is the number of crossings it replaces at once, $|m|+2$; see
Remark~\ref{rem:threestrands}.

\subsection{The family}

Fix $n\ge4$ and $2\le j\le n-2$, and write
\[
  a=\sigma_{j-1},\qquad b=\sigma_j,\qquad c=\sigma_{j+1}.
\]
For $\varepsilon\in\{1,-1\}$ and nonzero integers $p,q$, put
\[
  u_\varepsilon=a^{\varepsilon}c^{-\varepsilon},\qquad
  Y=u_\varepsilon\, b^{q}\, u_\varepsilon^{-1},\qquad
  w_\varepsilon(p,q)=b^{-p}\,u_\varepsilon\, b^{q}\,u_\varepsilon^{-1}\,
                    b^{p}\,u_\varepsilon\, b^{-q}\,u_\varepsilon^{-1}
                  =b^{-p}\,Y\,b^{p}\,Y^{-1},
\]
where $b^{q}$ denotes the word consisting of $|q|$ letters $b^{\pm1}$, and
$u_\varepsilon^{-1}=c^{\varepsilon}a^{-\varepsilon}$. The word $w_\varepsilon(p,q)$ has
length $2|p|+2|q|+8$ and uses only the three generators $a,b,c$. Since $a$ and $c$
commute, each two-letter block $u_\varepsilon^{\pm1}$ may be written in either order,
and we regard all $16$ resulting words as members of the family.

\begin{remark}
In particular, the hard braid
\[
  W=\sigma_2^{-2}\,\sigma_3^{-1}\sigma_1\,\sigma_2^{2}\,\sigma_1^{-1}\sigma_3\,
    \sigma_2^{2}\,\sigma_3^{-1}\sigma_1\,\sigma_2^{-2}\,\sigma_1^{-1}\sigma_3\in B_4
\]
is $w_{1}(2,2)$ with $j=2$, each block $u_1^{\pm1}$ being written with its negative
letter first. It is a shortest member of the family; Section~\ref{sec:diagram} shows it as
a diagram.
\end{remark}

\begin{remark}[Geometric meaning]\label{rem:geom}
It helps to keep in mind what these words do to the four strands in positions
$j-1,\dots,j+2$. Take $\varepsilon=1$. The word $b^{\mp p}$ twists the two inner strands
around each other. In $u_1=a\,c^{-1}$ the two outer strands move inwards, each passing
\emph{over} the inner strand next to it, so that they become adjacent; then $b^{q}$ twists
them around each other, and $u_1^{-1}$ takes them back. So $Y$ is a twist of the two outer
strands, performed above the two inner strands, and $w_1(p,q)$ is the commutator of a twist
of the inner pair and a twist of the outer pair. Since the outer pair lies entirely above
the inner pair, the two pairs can be pulled apart, and the two twists commute: this is the
reason why $w_1(p,q)$ is trivial, made precise in Lemma~\ref{lem:layer} below. For
$\varepsilon=-1$ the outer strands pass under the inner ones instead. The picture also
explains the opposite signs in $u_\varepsilon=a^{\varepsilon}c^{-\varepsilon}$: in
$a^{\varepsilon}c^{\varepsilon}$ one outer strand would pass over the inner pair and the
other under it. In this section we nevertheless prove triviality by explicit moves
(Lemma~\ref{lem:slide}), because these moves are needed in Section~\ref{sec:power}.
\end{remark}

\subsection{The slide lemma}

\begin{lemma}[slide]\label{lem:slide}
For every $\varepsilon\in\{1,-1\}$ and every nonzero integer $q$, the words $b^{\pm1}Y$
and $Yb^{\pm1}$ are related by a sequence of length-preserving moves: two power braid
moves (which are elementary braid moves when $|q|=1$), two elementary braid moves, and
$|q|$ commutations. In particular $b$ commutes with $Y$ in $B_n$.
\end{lemma}

\begin{proof}
Take first $\varepsilon=1$, so that $Y=a\,c^{-1}b^{q}c\,a^{-1}$, and the letter $b$.
Then
\[
\begin{aligned}
  b\,Y=b\,a\,(c^{-1}b^{q}c)\,a^{-1}
  &\ \longrightarrow\ b\,a\,(b\,c^{q}\,b^{-1})\,a^{-1}
  &&\text{power braid move,}\\
  &\ \longrightarrow\ (a\,b\,a)\,c^{q}\,b^{-1}a^{-1}
  &&\text{braid move } bab\to aba,\\
  &\ \longrightarrow\ a\,b\,c^{q}\,(a\,b^{-1}a^{-1})
  &&\text{$|q|$ commutations of $a$ with $c$,}\\
  &\ \longrightarrow\ a\,b\,c^{q}\,(b^{-1}a^{-1}b)
  &&\text{braid move } ab^{-1}a^{-1}\to b^{-1}a^{-1}b,\\
  &\ \longrightarrow\ a\,(c^{-1}b^{q}c)\,a^{-1}\,b=Y\,b
  &&\text{power braid move.}
\end{aligned}
\]
The braid moves have exponent patterns $(+,+,+)$ and $(+,-,-)$, which do not
alternate, and every step preserves length.

The remaining cases follow from this one by two symmetries, each of which maps legal moves
to legal moves and preserves lengths.
\begin{itemize}
  \item \emph{Exchanging $a$ and $c$}, and fixing $b$. This preserves which generators
    commute and which are adjacent. It carries $u_1=ac^{-1}$ to $ca^{-1}$, which is
    $u_{-1}=a^{-1}c$ written in the other order; so it carries the slide of $b$ across $Y$
    for $(\varepsilon,q)=(1,q)$ to the slide of $b$ across $Y$ for $(-1,q)$.
  \item \emph{Inverting every letter.} Reversing all exponents preserves the condition
    $x=z\ne y$, and transforms power braid moves into power braid moves. It carries the
    slide of $b$ across $Y$ for $(\varepsilon,q)$ to the slide of $b^{-1}$ across $Y$ for
    $(-\varepsilon,-q)$.
\end{itemize}
Since $q$ is an arbitrary nonzero integer, this covers both letters $b^{\pm1}$ and both
signs $\varepsilon$. Finally, every move is reversible, so the sequences also run from
$Yb^{\pm1}$ to $b^{\pm1}Y$.
\end{proof}

\subsection{Hardness}

\begin{theorem}\label{thm:hard}
Let $n\ge4$, $2\le j\le n-2$, $\varepsilon\in\{1,-1\}$, and let $p,q$ be integers with
$|p|\ge2$ and $|q|\ge2$. Then $w_\varepsilon(p,q)$ is a trivial braid word of length
$2|p|+2|q|+8$, and it is hard.
\end{theorem}

\begin{proof}
\emph{Triviality.} By Lemma~\ref{lem:slide}, $b$ commutes with $Y$, hence so does
$b^{p}$, and $w_\varepsilon(p,q)=b^{-p}\,Y\,b^{p}\,Y^{-1}=1$.

\emph{Hardness.} The idea is to exhibit a set of words containing $w_\varepsilon(p,q)$
which no move can leave and in which no deletion is possible. Let $\Omega$ be the set of cyclic words obtained from
$w_\varepsilon(p,q)$ by choosing the order of the two letters in each of its four blocks
$u_\varepsilon^{\pm1}$, together with their reversed inverses. Every word in $\Omega$
has the following shape, read cyclically: four maximal runs of the letter $b$ or
$b^{-1}$, of lengths $|p|,|q|,|p|,|q|$, each at least $2$, separated by four blocks of
exactly two letters, one a power of $a$ and one a power of $c$. We show that
$\Omega$ is closed under elementary moves and that no word of $\Omega$ admits a
deletion. Since deletions are the only moves that change the length, the empty word
is then unreachable.

\begin{itemize}
  \item \emph{No deletion.} Inside a run all letters coincide; inside a block the two
    letters belong to different generators; and at the junction of a run and a block
    the letters belong to $b$ and to $a$ or $c$. So no two adjacent letters are
    inverse to each other.
  \item \emph{Commutations.} The only commuting pair among $a,b,c$ is $\{a,c\}$, and
    letters of $a$ and $c$ are adjacent only inside a block. Exchanging them reorders
    a block and stays in $\Omega$.
  \item \emph{No braid move.} A braid move needs a subword $x\,y\,x'$ in which $x$
    and $x'$ belong to the same generator and $y$ to an adjacent one. Two letters of
    the same generator separated by exactly one letter do not occur: two letters $b$
    in different runs are separated by a block of two letters, two letters $b$ in
    the same run are separated only by letters $b$, and two letters of $a$, or of
    $c$, lie in different blocks and are separated by a run of length at least $2$.
\end{itemize}
Cyclic permutation and passage to the reversed inverse preserve the shape, so they
too stay within $\Omega$.
\end{proof}

\begin{remark}
The hypotheses are sharp in two respects.
\begin{enumerate}
  \item \emph{Signs.} If the two letters of $u$ have equal signs,
    $u=a^{\varepsilon}c^{\varepsilon}$, the word is not trivial: by machine, for all
    $p,q\in\{\pm1,\dots,\pm5,6\}$ in $B_4$, the word is trivial in all $242$ cases with opposite
    signs and in none of the $242$ cases with equal signs.
  \item \emph{Powers.} The hardness argument needs every run of $b$ to have length at
    least $2$. If $|p|=1$ or $|q|=1$, a braid move becomes available, and in every case
    we computed, with $p,q\in\{1,2,-2,3\}$ in $B_5$, the word is then untangled by elementary
    moves.
\end{enumerate}
The second statement is an observation, not a theorem.
\end{remark}

\begin{remark}
All words of the family lie in the subgroup generated by $\sigma_{j-1},\sigma_j,
\sigma_{j+1}$, which is isomorphic to $B_4$. The phenomenon is therefore a
four-strand phenomenon, embedded in every $B_n$ with $n\ge4$. Theorem~\ref{thm:hard}
gives infinitely many hard braids, the shortest of length $16$. On four strands there is no
shorter hard braid, by the exhaustive search described in Section~\ref{sec:verification};
on five strands there is none of length at most $12$.
\end{remark}

\section{A diagram of the hard braid}\label{sec:diagram}

Theorem~\ref{thm:hard} proves hardness algebraically. It is instructive to see the
same braid geometrically; this section illustrates our first contribution and prepares
the diagrammatic language used in Section~\ref{sec:three}. Every trivial braid word is the boundary word of a van
Kampen diagram over the presentation of $B_n$: a planar disc tiled by squares, with
boundary $\sigma_i\sigma_k\sigma_i^{-1}\sigma_k^{-1}$ for $|i-k|\ge2$, and hexagons,
with boundary $\sigma_i\sigma_k\sigma_i\sigma_k^{-1}\sigma_i^{-1}\sigma_k^{-1}$ for
$|i-k|=1$.

\emph{Labels.} We label each vertex by the permutation of the strands that it
represents, listing the strands from left to right, so that every edge interchanges
two strands in adjacent positions. We label the edge by these two strands and orient
it towards its endpoint of greater height, defined as follows.

\emph{Heights.} The number of positive
letters of a braid word minus the number of its negative letters is known in braid
theory as the \emph{exponent sum} of the braid, or its \emph{(signed) crossing number},
or its \emph{writhe}; see \cite{LisitsaSallesVernitski}, where it reappears as an
invariant learned by a neural network. It is unchanged by the defining relations, so it
is a homomorphism $B_n\to\mathbb Z$; it is also unchanged by cyclic permutations, by the
elementary moves and by the power braid moves, while passing to the reversed inverse
changes its sign. In particular, every trivial word has exponent sum $0$. The \emph{height}
of a vertex is the exponent sum of the word read along any path to it from a fixed base
vertex; it does not depend on the path, because the boundary word of every face is a
relator. We keep the geometric name because in Section~\ref{sec:three} we draw each
vertex at its height.

\emph{Reading words; tracks.} Walking along an edge whose
strands occupy positions $i$ and $i+1$ then reads $\sigma_i$ if the walk follows the
arrow and $\sigma_i^{-1}$ if it goes against it. Edges on opposite sides of a face
interchange the same two strands, and the chains of such edges form curves dual to the
edges, which we call \emph{tracks}.

\begin{figure}[ht]
\centering
\begin{tikzpicture}[scale=1.55,vtx/.style={font=\scriptsize,inner sep=1.2pt,fill=white,rounded corners=1pt},e/.style={postaction={decorate},decoration={markings,mark=at position 0.72 with {\arrow{Stealth[length=1.5mm]}}}},lb/.style={pos=0.32,fill=white,inner sep=0.4pt,font=\tiny}]
\fill[gray!12] (-2.00,0.00) -- (-1.00,0.00) -- (-1.00,-1.00) -- (-2.55,-2.55) -- (-3.57,-1.50) -- (-3.50,0.00) -- cycle;
\fill[gray!12] (-1.00,-1.00) -- (0.00,-1.00) -- (0.00,-2.00) -- (0.00,-3.50) -- (-1.50,-3.57) -- (-2.55,-2.55) -- cycle;
\fill[orange!22] (0.00,-2.00) -- (0.00,-1.00) -- (1.00,-1.00) -- (2.55,-2.55) -- (1.50,-3.57) -- (0.00,-3.50) -- cycle;
\fill[gray!12] (1.00,-1.00) -- (1.00,0.00) -- (2.00,0.00) -- (3.50,0.00) -- (3.57,-1.50) -- (2.55,-2.55) -- cycle;
\fill[orange!22] (2.00,0.00) -- (1.00,0.00) -- (1.00,1.00) -- (2.55,2.55) -- (3.57,1.50) -- (3.50,0.00) -- cycle;
\fill[orange!22] (1.00,1.00) -- (0.00,1.00) -- (0.00,2.00) -- (0.00,3.50) -- (1.50,3.57) -- (2.55,2.55) -- cycle;
\fill[gray!12] (0.00,2.00) -- (0.00,1.00) -- (-1.00,1.00) -- (-2.55,2.55) -- (-1.50,3.57) -- (0.00,3.50) -- cycle;
\fill[orange!22] (-1.00,1.00) -- (-1.00,0.00) -- (-2.00,0.00) -- (-3.50,0.00) -- (-3.57,1.50) -- (-2.55,2.55) -- cycle;
\fill[teal!18] (-1.00,1.00) -- (-1.00,0.00) -- (0.00,0.00) -- (0.00,1.00) -- cycle;
\fill[teal!18] (-1.00,0.00) -- (-1.00,-1.00) -- (0.00,-1.00) -- (0.00,0.00) -- cycle;
\fill[teal!18] (0.00,-1.00) -- (1.00,-1.00) -- (1.00,0.00) -- (0.00,0.00) -- cycle;
\fill[teal!18] (1.00,0.00) -- (1.00,1.00) -- (0.00,1.00) -- (0.00,0.00) -- cycle;
\draw[red!70!black,line width=1.1pt,dash pattern=on 4pt off 2.5pt] (-2.75,0.00) -- (-1.77,-1.77) -- (0.00,-2.75) -- (1.77,-1.77) -- (2.75,0.00) -- (1.77,1.77) -- (0.00,2.75) -- (-1.77,1.77) -- cycle;
\node[vtx] (v15) at (-1.00,1.00) {1423};
\node[vtx] (v0) at (-1.00,0.00) {1432};
\node[vtx] (v3) at (-1.00,-1.00) {1423};
\node[vtx] (v4) at (0.00,-1.00) {4123};
\node[vtx] (v7) at (1.00,-1.00) {1423};
\node[vtx] (v8) at (1.00,0.00) {1432};
\node[vtx] (v11) at (1.00,1.00) {1423};
\node[vtx] (v12) at (0.00,1.00) {4123};
\node[vtx] (v1) at (-2.00,0.00) {1342};
\node[vtx] (v5) at (0.00,-2.00) {4213};
\node[vtx] (v9) at (2.00,0.00) {1342};
\node[vtx] (v13) at (0.00,2.00) {4213};
\node[vtx] (ctr) at (0.00,0.00) {4132};
\node[vtx] (uv1) at (-3.50,0.00) {1324};
\node[vtx] (uv3) at (-2.55,-2.55) {1243};
\node[vtx] (t1) at (-3.57,-1.50) {1234};
\node[vtx] (uv5) at (0.00,-3.50) {2413};
\node[vtx] (t3) at (-1.50,-3.57) {2143};
\node[vtx] (uv7) at (2.55,-2.55) {1243};
\node[vtx] (t5) at (1.50,-3.57) {2143};
\node[vtx] (uv9) at (3.50,0.00) {1324};
\node[vtx] (t7) at (3.57,-1.50) {1234};
\node[vtx] (uv11) at (2.55,2.55) {1243};
\node[vtx] (t9) at (3.57,1.50) {1234};
\node[vtx] (uv13) at (0.00,3.50) {2413};
\node[vtx] (t11) at (1.50,3.57) {2143};
\node[vtx] (uv15) at (-2.55,2.55) {1243};
\node[vtx] (t13) at (-1.50,3.57) {2143};
\node[vtx] (t15) at (-3.57,1.50) {1234};
\draw[e] (v0) -- (ctr) node[lb] {14};
\draw[e] (ctr) -- (v12) node[lb] {23};
\draw[e] (v4) -- (ctr) node[lb] {23};
\draw[e] (ctr) -- (v8) node[lb] {14};
\draw[e] (t1) -- (uv1) node[lb] {23};
\draw[e] (uv3) -- (t1) node[lb] {34};
\draw[e] (uv11) -- (t11) node[lb] {12};
\draw[e] (uv13) -- (t11) node[lb] {14};
\draw[e] (t13) -- (uv13) node[lb] {14};
\draw[e] (uv15) -- (t13) node[lb] {12};
\draw[e] (uv1) -- (t15) node[lb] {23};
\draw[e] (uv15) -- (t15) node[lb] {34};
\draw[e] (uv3) -- (t3) node[lb] {12};
\draw[e] (t3) -- (uv5) node[lb] {14};
\draw[e] (uv5) -- (t5) node[lb] {14};
\draw[e] (uv7) -- (t5) node[lb] {12};
\draw[e] (uv7) -- (t7) node[lb] {34};
\draw[e] (t7) -- (uv9) node[lb] {23};
\draw[e] (uv11) -- (t9) node[lb] {34};
\draw[e] (uv9) -- (t9) node[lb] {23};
\draw[e] (v1) -- (uv1) node[lb] {24};
\draw[e] (v11) -- (uv11) node[lb] {24};
\draw[e] (v13) -- (uv13) node[lb] {24};
\draw[e] (v15) -- (uv15) node[lb] {24};
\draw[e] (v3) -- (uv3) node[lb] {24};
\draw[e] (v5) -- (uv5) node[lb] {24};
\draw[e] (v7) -- (uv7) node[lb] {24};
\draw[e] (v9) -- (uv9) node[lb] {24};
\draw[e] (v0) -- (v1) node[lb] {34};
\draw[e] (v0) -- (v15) node[lb] {23};
\draw[e] (v3) -- (v0) node[lb] {23};
\draw[e] (v12) -- (v11) node[lb] {14};
\draw[e] (v8) -- (v11) node[lb] {23};
\draw[e] (v12) -- (v13) node[lb] {12};
\draw[e] (v15) -- (v12) node[lb] {14};
\draw[e] (v3) -- (v4) node[lb] {14};
\draw[e] (v4) -- (v5) node[lb] {12};
\draw[e] (v4) -- (v7) node[lb] {14};
\draw[e] (v7) -- (v8) node[lb] {23};
\draw[e] (v8) -- (v9) node[lb] {34};
\node[font=\small] at (-4.05,1.78) {$\ast$};
\node[font=\small] at (-1.27,0.21) {$\dagger$};
\end{tikzpicture}
\caption{A diagram whose boundary word is the hard braid $W$. Each vertex is
labelled by a permutation of the strands $1,2,3,4$, listed from left to right; thus
$1432$ has strand $1$ in the first position. Each edge is labelled by the two strands it
interchanges and points towards greater height. The squares (teal) form a $2\times2$
grid; the hexagons form a ring; the four orange hexagons are those whose
inner side has its lowest vertex in the middle, see the text. The dashed curve is the closed
track of strands $2$ and $4$. Reading the outer boundary counterclockwise from the
vertex $1234$ marked $\ast$, and recording at each edge $\sigma_i$ or $\sigma_i^{-1}$
according to the positions $i,i+1$ of the two strands and the direction of the arrow,
gives $W$. Reading the boundary of the grid with its
spikes counterclockwise from the vertex $1432$ marked $\dagger$ gives the inner word $r$.}
\label{fig:basketW}
\end{figure}
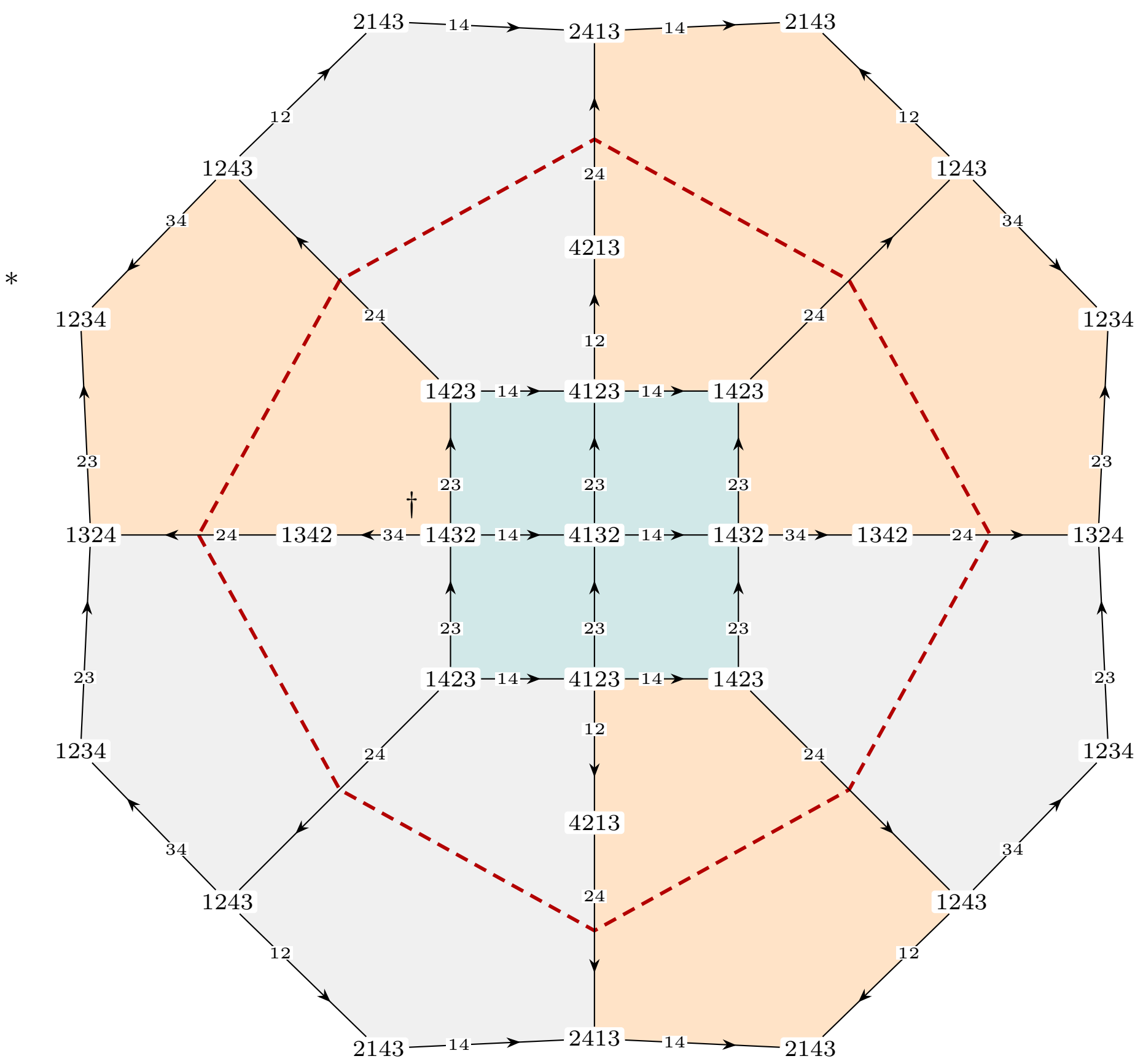

Figure~\ref{fig:basketW} shows a diagram whose boundary word is $W$. Its outer boundary
is drawn as an octagon whose sides each consist of two edges: the four horizontal and
vertical sides carry the four runs $\sigma_2^{\pm2}$ of $W$, and the four diagonal
sides carry the four blocks $u^{\pm1}$. Inside, a ring of eight hexagons surrounds a
$2\times2$ grid of four squares, with a spike from the middle of each side of the grid.
The edges of the grid interchange strands $1$ and $4$ or strands $2$ and $3$, so the
grid is crossed by two tracks of each kind, and each track of one kind crosses both
tracks of the other. A closed track, interchanging strands $2$ and $4$, runs round the
ring of hexagons.

If a word can be untangled without increasing its length, running the untangling
backwards builds a diagram for it in which every track is an arc with both ends on the
boundary. A closed track therefore indicates that the diagram cannot have been built
that way. This alone does not prove that a word is hard, since a word has many
diagrams, and in most diagrams with a closed track the boundary word can be realised by
another diagram without one; but $W$ was first found by searching for diagrams whose
closed tracks cannot be removed, and Figure~\ref{fig:basketW} is the diagram in which it
was found. Theorem~\ref{thm:hard} explains why such a diagram exists.

\paragraph{The inner word.}
The hexagons of the ring play two different roles, shown in grey and orange in
Figure~\ref{fig:basketW}, and the clearest way to see the difference is to compare $W$
with the word read along the inner side of the ring. The grid of squares together with
its four spikes forms a region $R$, and the boundary of $R$ is the inner boundary of the
ring. Walk once round the boundary of $R$ counterclockwise, starting at the vertex $1432$
in the middle of the left side of the grid, marked $\dagger$; whenever the walk reaches a
spike it goes out along it to its free end and comes back along the same edge. The walk
has sixteen steps, and recording a letter at each step by the rule above gives the
\emph{inner word}
\[
r=\sigma_2\sigma_2^{-1}\cdot\sigma_3^{-1}\sigma_1\cdot\sigma_2\sigma_2^{-1}\cdot
\sigma_1\sigma_3\cdot\sigma_2\sigma_2^{-1}\cdot\sigma_3\sigma_1^{-1}\cdot
\sigma_2\sigma_2^{-1}\cdot\sigma_1^{-1}\sigma_3^{-1},
\]
where the dots separate the four excursions along spikes, each reading
$\sigma_2\sigma_2^{-1}$, from the four turns round the corners of the grid, each reading
a letter $\sigma_1^{\pm1}$ and a letter $\sigma_3^{\pm1}$.

The word $r$ is easy to untangle. Deleting the four pairs $\sigma_2\sigma_2^{-1}$ leaves
$\sigma_3^{-1}\sigma_1^{2}\sigma_3^{2}\sigma_1^{-2}\sigma_3^{-1}$, the boundary word of
the grid. Since $\sigma_1$ and $\sigma_3$ commute, four commutations turn
$\sigma_1^{2}\sigma_3^{2}$ into $\sigma_3^{2}\sigma_1^{2}$, after which every letter
cancels against its neighbour. Thus $r$ is carried to the empty word without its length
ever increasing, in accordance with the fact that $R$ contains no closed track.

Now compare $r$ with $W$. Each hexagon of the ring meets the boundary of $R$ in two
consecutive edges of the walk from $\dagger$ and the outer boundary in two consecutive
edges of the walk from $\ast$; the hexagon containing both $\dagger$ and $\ast$ supplies
the first and the last letter of each word, and the other seven supply letters $2$
and~$3$, $4$ and~$5$, and so on up to $14$ and~$15$. Going round a hexagon from the first vertex
of its inner side to the last vertex of its outer side, one route follows the inner side and then an edge
crossing the closed track, and the other route follows an edge crossing the closed track
and then the outer side. The two routes read $\sigma_i\sigma_k\sigma_i$ and
$\sigma_k\sigma_i\sigma_k$ up to signs, with $|i-k|=1$, so the inner and outer sides
each read $\sigma_i$ then $\sigma_k$, up to signs. Hence the $m$-th letters of $r$ and of
$W$ are powers of the same generator, for every $m$; the two words can differ only in
their signs.

The signs are governed by heights. In a grey hexagon the inner and outer sides both rise
or both fall: for instance the hexagon below the grid on the left has heights $-1,0,1$
along its inner side and $0,1,2$ along its outer side. Its two letters are then the same
in $r$ and in $W$. In an orange hexagon the inner side falls and then rises, so that its
lowest vertex is in the middle, while the outer side rises and then falls: heights $1,0,1$
against $2,3,2$, or $3,2,3$ against $4,5,4$. Both of its letters then have opposite signs
in $r$ and in $W$. Consequently $W$ is obtained from $r$ by reversing the signs of the
eight letters that lie in orange hexagons, shaded below:
\[
\setlength{\arraycolsep}{1.6pt}
\begin{array}{r*{16}{c}}
 & \scriptstyle1&\scriptstyle2&\scriptstyle3&\scriptstyle4&\scriptstyle5&\scriptstyle6&\scriptstyle7&\scriptstyle8&\scriptstyle9&\scriptstyle10&\scriptstyle11&\scriptstyle12&\scriptstyle13&\scriptstyle14&\scriptstyle15&\scriptstyle16\\[2pt]
r\colon & \fl{\sigma_2} & \sigma_2^{-1} & \sigma_3^{-1} & \sigma_1 & \sigma_2 & \fl{\sigma_2^{-1}} & \fl{\sigma_1} & \sigma_3 & \sigma_2 & \fl{\sigma_2^{-1}} & \fl{\sigma_3} & \fl{\sigma_1^{-1}} & \fl{\sigma_2} & \sigma_2^{-1} & \sigma_1^{-1} & \fl{\sigma_3^{-1}}\\[3pt]
W\colon & \fl{\sigma_2^{-1}} & \sigma_2^{-1} & \sigma_3^{-1} & \sigma_1 & \sigma_2 & \fl{\sigma_2} & \fl{\sigma_1^{-1}} & \sigma_3 & \sigma_2 & \fl{\sigma_2} & \fl{\sigma_3^{-1}} & \fl{\sigma_1} & \fl{\sigma_2^{-1}} & \sigma_2^{-1} & \sigma_1^{-1} & \fl{\sigma_3}
\end{array}
\]
Each excursion $\sigma_2\sigma_2^{-1}$ along a spike has exactly one letter reversed and
becomes one of the four runs $\sigma_2^{\pm2}$ of $W$, and three of the four corner turns
have one or both letters reversed.

These eight sign changes are the whole difference between an easy word and a hard one:
$r$ untangles by elementary moves, while $W$ does not, by Theorem~\ref{thm:hard}. If every
hexagon of the ring were grey, the outer word would coincide with $r$. So the ring of
hexagons threaded by the closed track carries the inner word to the outer one, and the
orange hexagons are the places where it changes it. The heights of individual vertices
play no role in the proofs of Section~\ref{sec:family}; in the figure they matter only
through the arrows, which record them, and through this comparison. We stress that the
comparison describes this particular diagram; we do not claim here a general statement
relating hard words to inner words of rings of hexagons.

\section{Three strands}\label{sec:three}

This section establishes our second contribution: hard braids need at least four
strands, so the family of Section~\ref{sec:family} is minimal in the number of strands.
On two strands this is immediate, since a nonempty word in $\sigma_1^{\pm1}$ with exponent
sum $0$ contains two adjacent letters of opposite signs, which can be deleted. The case
of three strands is the content of the following theorem.

\begin{theorem}\label{thm:three}
There are no hard braids in $B_3$. That is, every trivial word in
$\sigma_1^{\pm1},\sigma_2^{\pm1}$ can be transformed into the empty word by
deletions $\sigma_i^{\varepsilon}\sigma_i^{-\varepsilon}$ and braid moves with
non-alternating exponents, together with cyclic permutations.
\end{theorem}

\paragraph{Strategy.}
The proof is by induction on the area of the word, the least number of relators needed to
express it. A cyclically reduced trivial word bounds a van Kampen diagram with the least
possible number of faces. We show that such a diagram always has a face three of whose
consecutive edges lie on the boundary of the diagram, with nothing attached to the two
vertices between them (an \emph{exposed half-face}, Lemma~\ref{lem:exposed}). The braid move
that replaces these three edges by the other three edges of the face removes the face, and
so lowers the area without changing the length. The existence of an exposed half-face is
a counting argument (Lemma~\ref{lem:euler}) based on one local property of minimal diagrams
on three strands (Lemma~\ref{lem:corners}), which is exactly the property that fails on four
strands.

Throughout this section we write $a=\sigma_1$, $b=\sigma_2$ and
$R=aba\,b^{-1}a^{-1}b^{-1}$, so that $B_3=\langle a,b\mid R\rangle$, and we use
uppercase letters for inverses when convenient.

\paragraph{Area and diagrams.}
We recall standard facts, for which we refer to \cite[Ch.~V]{LyndonSchupp}. The
\emph{area} $A(w)$ of a trivial word $w$ is the least $N$ such that $w$ is equal, in
the free group on $a,b$, to a product of $N$ conjugates of $R^{\pm1}$. Since cyclic
permutation is conjugation, and free reduction does not change the element of the
free group, neither operation changes the area, and a word of area $0$ is freely
trivial. A \emph{van Kampen diagram} for $w$ is a finite, connected, simply connected
$2$-complex $D$ in the plane whose edges are oriented and labelled $a$ or $b$, such
that the boundary of every face reads a cyclic permutation of $R^{\pm1}$ and the
boundary cycle of $D$, read from a base vertex, reads $w$. By van Kampen's lemma, $A(w)$
is the least number of faces of a diagram for $w$. A diagram with $A(w)$ faces is
\emph{reduced}: it contains no two distinct faces $\Pi_1,\Pi_2$ with a common edge $e$
such that the boundary cycles of $\Pi_1$ and $\Pi_2$, each read starting along $e$,
have the same label, since such a pair can be removed, lowering the number of faces.

Let $D$ be a diagram for a trivial word $w$. Since $D$ is simply connected, there is a
map $\phi$ from the vertices of $D$ to $B_3$, taking the base vertex to $1$, such that
$\phi(y)=\phi(x)\sigma$ whenever an edge labelled $\sigma$ runs from $x$ to $y$. We
write $\pi(x)\in\mathfrak S_3$ for the image of $\phi(x)$ and $h(x)\in\mathbb Z$ for
its exponent sum, or signed crossing number, which is the \emph{height} of $x$ as in
Section~\ref{sec:diagram}. Every edge points from a vertex to a vertex one higher. At a vertex $x$, an edge \emph{points in} if it ends at $x$ and \emph{points
out} if it starts at $x$.

\begin{lemma}\label{lem:faces}
Every face $\Pi$ of $D$ is an embedded hexagon. Its vertices have heights
$m,m+1,m+2,m+3,m+2,m+1$ in cyclic order, for some $m$. The vertex of height $m$, the
\emph{source} of $\Pi$, is joined to the vertex of height $m+3$, the \emph{sink}, by two
increasing paths along $\partial\Pi$, labelled $aba$ and $bab$; we call them the
\emph{sides} of $\Pi$.
\end{lemma}

\begin{proof}
Reading $\partial\Pi$ from a suitable vertex $x$ in a suitable direction gives $R$, so
the vertices of $\Pi$ are mapped by $\phi$ to $g$, $ga$, $gab$, $gaba$, $gabab^{-1}$ and
$gabab^{-1}a^{-1}$, where $g=\phi(x)$. Their images under $\pi$ are
$\pi(x)$ multiplied on the right by $1$, $s_1$, $s_1s_2$, $s_1s_2s_1$, $s_2s_1$ and
$s_2$, where $s_1,s_2$ are the adjacent transpositions; these six permutations are
distinct, so the six vertices of $\Pi$ are distinct. The heights are read off from the
exponent sums of the prefixes of $R$, which are $0,1,2,3,2,1$. The path
$g\to ga\to gab\to gaba$ is labelled $aba$, and the path from $g$ through
$gabab^{-1}a^{-1}=gb$ and $gabab^{-1}=gba$ to $gaba=gbab$ is labelled $bab$.
\end{proof}

\begin{figure}[ht]
\centering
\begin{minipage}[b]{0.40\textwidth}\centering
\begin{tikzpicture}[scale=1.12,vtx/.style={font=\scriptsize,inner sep=1.2pt,fill=white,rounded corners=1pt},e/.style={postaction={decorate},decoration={markings,mark=at position 0.62 with {\arrow{Stealth[length=1.5mm]}}}},lb/.style={pos=0.4,fill=white,inner sep=0.4pt,font=\tiny}]
\fill[gray!12] (1.50,0.00) -- (3.00,1.05) -- (3.00,2.10) -- (1.50,3.15) -- (0.00,2.10) -- (0.00,1.05) -- cycle;

\node[vtx] (s) at (1.50,0.00) {123};
\node[vtx] (r1) at (3.00,1.05) {213};
\node[vtx] (r2) at (3.00,2.10) {231};
\node[vtx] (k) at (1.50,3.15) {321};
\node[vtx] (l1) at (0.00,1.05) {132};
\node[vtx] (l2) at (0.00,2.10) {312};
\draw[e] (s) -- (r1) node[lb] {12};
\draw[e] (r1) -- (r2) node[lb] {13};
\draw[e] (r2) -- (k) node[lb] {23};
\draw[e] (s) -- (l1) node[lb] {23};
\draw[e] (l1) -- (l2) node[lb] {13};
\draw[e] (l2) -- (k) node[lb] {12};
\node[font=\footnotesize,below=1pt] at (s.south) {source};\node[font=\footnotesize,above=1pt] at (k.north) {sink};\node[font=\footnotesize] at (3.55,1.6) {$aba$};\node[font=\footnotesize] at (-0.55,1.6) {$bab$};
\end{tikzpicture}\\[3pt](a)
\end{minipage}\hfill
\begin{minipage}[b]{0.56\textwidth}\centering
\begin{tikzpicture}[scale=1.12,vtx/.style={font=\scriptsize,inner sep=1.2pt,fill=white,rounded corners=1pt},e/.style={postaction={decorate},decoration={markings,mark=at position 0.62 with {\arrow{Stealth[length=1.5mm]}}}},lb/.style={pos=0.4,fill=white,inner sep=0.4pt,font=\tiny}]
\fill[red!10] (0.00,3.15) -- (0.00,2.10) -- (-0.56,1.05) -- (-1.38,0.00) -- (-2.31,1.05) -- (-2.31,2.10) -- cycle;
\fill[red!10] (0.00,3.15) -- (0.00,2.10) -- (0.56,1.05) -- (1.38,0.00) -- (2.31,1.05) -- (2.31,2.10) -- cycle;

\node[vtx] (x) at (0.00,3.15) {321};
\node[vtx] (w) at (0.00,2.10) {231};
\node[vtx] (La) at (-0.56,1.05) {213};
\node[vtx] (Lg) at (-1.38,0.00) {123};
\node[vtx] (Lb) at (-2.31,1.05) {132};
\node[vtx] (Lba) at (-2.31,2.10) {312};
\node[vtx] (Ra) at (0.56,1.05) {213};
\node[vtx] (Rg) at (1.38,0.00) {123};
\node[vtx] (Rb) at (2.31,1.05) {132};
\node[vtx] (Rba) at (2.31,2.10) {312};
\draw[e] (w) -- (x) node[lb] {23};
\draw[e] (La) -- (w) node[lb] {13};
\draw[e] (Lg) -- (La) node[lb] {12};
\draw[e] (Lg) -- (Lb) node[lb] {23};
\draw[e] (Lb) -- (Lba) node[lb] {13};
\draw[e] (Lba) -- (x) node[lb] {12};
\draw[e] (Ra) -- (w) node[lb] {13};
\draw[e] (Rg) -- (Ra) node[lb] {12};
\draw[e] (Rg) -- (Rb) node[lb] {23};
\draw[e] (Rb) -- (Rba) node[lb] {13};
\draw[e] (Rba) -- (x) node[lb] {12};
\node[font=\footnotesize,above=1pt] at (x.north) {$x$};\node[font=\footnotesize] at (0.24,2.86) {$e$};\node[font=\small] at (-1.05,1.55) {$\Pi_1$};\node[font=\small] at (1.05,1.55) {$\Pi_2$};
\end{tikzpicture}\\[3pt](b)
\end{minipage}
\caption{(a) A face of a diagram over $B_3$, with each vertex drawn at its height.
Vertices are labelled by permutations of the strands and edges by the two strands they
interchange, as in Figure~\ref{fig:basketW}. The two sides, from the source $123$ to the
sink $321$, read $aba$ and $bab$; the other four vertices carry side corners.
(b) Two faces with sink corners at the same vertex $x$, sharing the edge $e$. They are
mirror images of each other, with the same labels, and so they cancel; by
Lemma~\ref{lem:corners} this does not happen in a reduced diagram.}
\label{fig:corners}
\end{figure}
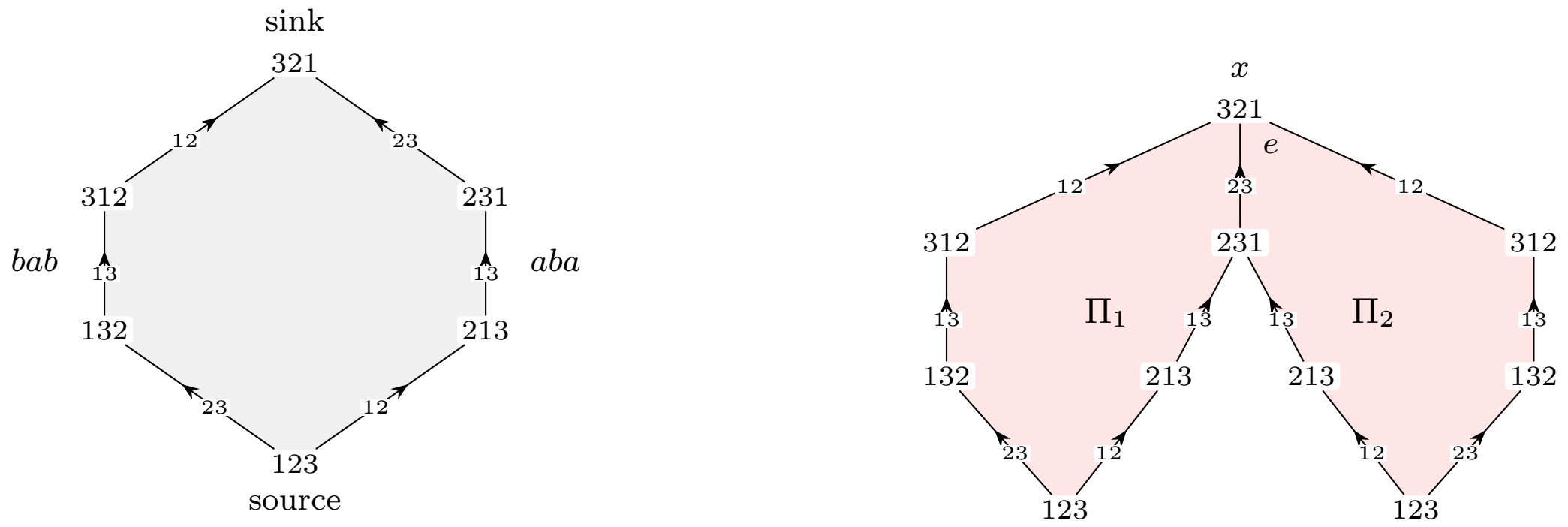

Since faces are embedded, a face containing a vertex $x$ meets $x$ in a single
\emph{corner}, bounded by the two edges of the face at $x$. The corner is a \emph{sink
corner} if $x$ is the sink of the face, so that both edges point in, a \emph{source
corner} if $x$ is its source, and a \emph{side corner} otherwise, when one edge points in
and one points out. By Lemma~\ref{lem:faces}, every face has one source corner, one sink
corner and four side corners.
Figure~\ref{fig:corners}(a) shows a face. As in Figure~\ref{fig:basketW}, vertices are
labelled by permutations of the strands and edges by the two strands they interchange,
but here each vertex is drawn at its height, so that every arrow points upwards.

\begin{lemma}\label{lem:corners}
If $D$ is reduced, two corners at the same vertex that share an edge are never both
sink corners, nor both source corners.
\end{lemma}

\begin{proof}
Two corners at a vertex $x$ sharing an edge $e$ belong to two distinct faces $\Pi_1$,
$\Pi_2$, because a face has only one corner at $x$. Suppose that both are sink corners,
and that $e$ is labelled $a$; the other cases are symmetric. By Lemma~\ref{lem:faces},
$e$ is the last edge of the side of $\Pi_i$ labelled $aba$, for $i=1,2$. Hence the
boundary cycle of $\Pi_i$, read starting along $e$, follows $e$ to the sink, descends
the side labelled $bab$ and climbs the first two edges of the other side, and so reads
$a\,b^{-1}a^{-1}b^{-1}\,ab$ for both $i=1$ and $i=2$. Thus $\Pi_1$ and $\Pi_2$ are mirror images
of each other across $e$, as in Figure~\ref{fig:corners}(b), and form a reducible pair,
which is impossible in a reduced diagram. Source corners are treated in
the same way, with $e$ the first edge of a side.
\end{proof}

\begin{remark}
Lemma~\ref{lem:corners} is the only place where we use that there are three strands.
The point is that the two edges at the sink of a face carry different letters, so in
$B_3$ a sink corner at a given vertex is determined by either of its edges. In $B_n$
with $n\ge4$, two faces with sink corners at $x$ sharing an edge labelled $\sigma_i$
may have their other edges labelled by different letters $\sigma_k\ne\sigma_l$, and then
they do not cancel. The remaining steps of the argument use only the shape of the faces,
and they adapt to diagrams with square faces, each of which also has one source corner, one
sink corner and increasing sides; an exposed side of a square is then a commutation move.
We do not need this adaptation and omit its details. Granting it, in every diagram for a
hard braid with the least number of faces, two faces with sink corners, or two with source
corners, share an edge at some vertex: this is the feature of four strands that hard braids
exploit. In the diagram of Figure~\ref{fig:basketW} this happens at five places; for
instance, at a corner $1423$ of the grid, a square of the grid and a hexagon of the ring
share an edge and both have their source there.
\end{remark}

Let $B$ be a \emph{disc diagram}: a diagram, or a part of one, homeomorphic to a closed
disc, so that its boundary is a simple closed path. For an interior vertex $x$ of $B$,
let $c(x)$ be the number of side corners at $x$, which is the number of changes between
edges pointing in and edges pointing out as one goes once round $x$. For a vertex $x$ on
$\partial B$, list the edges of $B$ at $x$ in the order in which they are met going
round $x$ inside $B$ from one boundary edge to the other; consecutive edges in this list
bound the corners of $B$ at $x$, one for each pair. Let $s(x)$ be the number of side
corners among them, which is the number of changes between in and out along the list.

\begin{lemma}\label{lem:euler}
If every face of the disc diagram $B$ is a hexagon with one source corner, one sink
corner and four side corners, then
\[
\sum_{x\ \mathrm{interior}}\bigl(2-c(x)\bigr)+\sum_{x\in\partial B}\bigl(1-s(x)\bigr)=2 .
\]
If moreover $B$ satisfies the conclusion of Lemma~\ref{lem:corners}, then $c(x)\ge2$
for every interior vertex $x$, and hence $\sum_{x\in\partial B}(1-s(x))\ge2$.
\end{lemma}

\begin{proof}
Let $B$ have $V$ vertices, of which $V_\partial$ lie on $\partial B$, $E$ edges, of
which $E_\partial$ lie on $\partial B$, and $F$ faces. Since $\partial B$ is a simple
closed path, $V_\partial=E_\partial$. Counting side corners face by face gives
$\sum c(x)+\sum s(x)=4F$, and counting pairs of a face and an edge on its boundary gives
$6F=2E-E_\partial$. Hence the left-hand side equals
\[
2(V-V_\partial)+V_\partial-4F=2V-E_\partial-4F=2V-2E+2F=2,
\]
by Euler's formula for the disc. For the second statement, $c(x)$ is even, since it
counts changes round a cycle. If $c(x)=0$, then all edges at $x$ point in, or all point
out, and all corners at $x$ are sink corners, or all are source corners. An interior
vertex has at least two edges, since otherwise the face containing its only corner would
contain this edge twice in its boundary, and so it has at least two corners, two of which
share an edge. This contradicts Lemma~\ref{lem:corners}; hence $c(x)\ge2$.
\end{proof}

Going round $\partial B$, the height changes by $\pm1$ at each step. A boundary vertex
whose two boundary edges both point in is a \emph{peak}, and one whose two boundary edges
both point out is a \emph{valley}; peaks and valleys are \emph{turning points}, and the
other boundary vertices are \emph{slopes}. Peaks and valleys alternate round
$\partial B$, and there is at least one of each, so $\partial B$ is the union of
\emph{segments}: maximal increasing arcs, each running from a valley up to the next peak.
Each slope is an inner vertex of exactly one segment, and each turning point is an end of
exactly two. Along the list of edges at a boundary vertex $x$, the number $s(x)$ of
changes is odd if $x$ is a slope and even if $x$ is a turning point. We call a turning
point with $s(x)=0$ a \emph{tip}. At a tip all corners are sink corners, or all are
source corners, so under the conclusion of Lemma~\ref{lem:corners} a tip has exactly one
corner: a tip has degree $2$ and is the sink, if it is a peak, or the source, if it is a
valley, of the only face containing it. A turning point that is not a tip has
$s(x)\ge2$.

An \emph{exposed half-face} of $B$ is a path of three consecutive edges in the boundary
of a face of $B$ which are also consecutive edges of $\partial B$, and whose two inner
vertices have degree $2$ in $B$.

\begin{lemma}\label{lem:exposed}
Let $B$ be a disc diagram with at least one face, satisfying the conclusions of
Lemmas~\ref{lem:faces} and~\ref{lem:corners}. Then $\partial B$ has two distinct
segments each of which contains an exposed half-face whose two inner vertices are slopes
of that segment.
\end{lemma}

\begin{proof}
We first show, by weighing the boundary vertices, that at least two segments are as simple
as possible (\emph{full}, in the terminology below), and then that every full segment
contains an exposed half-face.

For a segment $I$ with ends $p$ and $q$, let
\[
\tau(I)=\sum_{x\text{ a slope of }I}\bigl(1-s(x)\bigr)
+\tfrac12\bigl(1-s(p)\bigr)+\tfrac12\bigl(1-s(q)\bigr).
\]
Each boundary vertex is counted with total weight $1$, so $\sum_I\tau(I)=
\sum_{x\in\partial B}(1-s(x))\ge2$ by Lemma~\ref{lem:euler}. Each slope contributes at
most $0$, an end which is a tip contributes $\frac12$, and an end which is not a tip
contributes at most $-\frac12$. So $\tau(I)$ is an integer at most $1$, with
equality exactly when both ends of $I$ are tips and every slope of $I$ has
$s(x)=1$; we then call $I$ \emph{full}. Since the values sum to at least $2$, at
least two segments are full. It remains to show that every full segment contains an
exposed half-face whose inner vertices are slopes of the segment.

Let $I$ be full, and list its vertices downwards as $x_0,x_1,\dots,x_L$, where the
peak $x_0$ and the valley $x_L$ are tips. For $0\le j<L$ let $b_j$ be the boundary edge
from $x_{j+1}$ up to $x_j$, and $\Pi_j$ the unique face containing it. A boundary vertex
of degree $2$ has a single corner, formed by its two boundary edges, so if $x_{j+1}$ has
degree $2$ then $\Pi_j=\Pi_{j+1}$, and $b_j,b_{j+1}$ are consecutive edges of
$\partial\Pi_j$. The peak $x_0$ is a tip, so it is the sink of $\Pi_0$. We show:
\emph{if $x_j$ is the sink of $\Pi_j$, then either $\Pi_j$ has an exposed half-face whose
inner vertices are slopes of $I$, or $x_k$ is the sink of $\Pi_k$ for some
$k$ with $j<k<L$.} Since $L$ is finite, this proves the claim.

Suppose that $x_j$ is the sink of $\Pi_j$. Let $\ell\ge1$ be largest such that $j+\ell\le
L$ and $x_{j+1},\dots,x_{j+\ell-1}$ all have degree $2$. Then $b_j,\dots,b_{j+\ell-1}$
form an increasing path in $\partial\Pi_j$ which ends at the sink $x_j$. By
Lemma~\ref{lem:faces}, such a path lies on a side of $\Pi_j$; so $\ell\le3$, and $\ell=3$
only if the path starts at the source of $\Pi_j$. If $\ell=3$, then $b_j,b_{j+1},b_{j+2}$
is an exposed half-face of $\Pi_j$, and its inner vertices $x_{j+1},x_{j+2}$ are slopes
of $I$, since $j+2<j+3\le L$.

If $\ell\le2$, then $y=x_{j+\ell}$ is not the source of $\Pi_j$, so the corner of $\Pi_j$
at $y$ is a side corner. If $y=x_L$, this contradicts the fact that $x_L$ is a tip, whose
only corner is a source corner. So $y$ is a slope of $I$, and by the maximality of
$\ell$ it has degree at least $3$, hence at least two corners. The list of edges at $y$
begins with $b_{j+\ell-1}$, which points out of $y$, followed by the other edge of
$\Pi_j$ at $y$, which points in, and it ends with $b_{j+\ell}$, which points in. Since
$s(y)=1$, every edge in the list after the first points in, so every corner at $y$ other
than that of $\Pi_j$ is a sink corner. There is at least one such corner, and by
Lemma~\ref{lem:corners} there is at most one. Hence $y$ has exactly two corners, and the
last of them, which is the corner of $\Pi_{j+\ell}$, is a sink corner: $x_{j+\ell}$ is the
sink of $\Pi_{j+\ell}$, as required.
\end{proof}

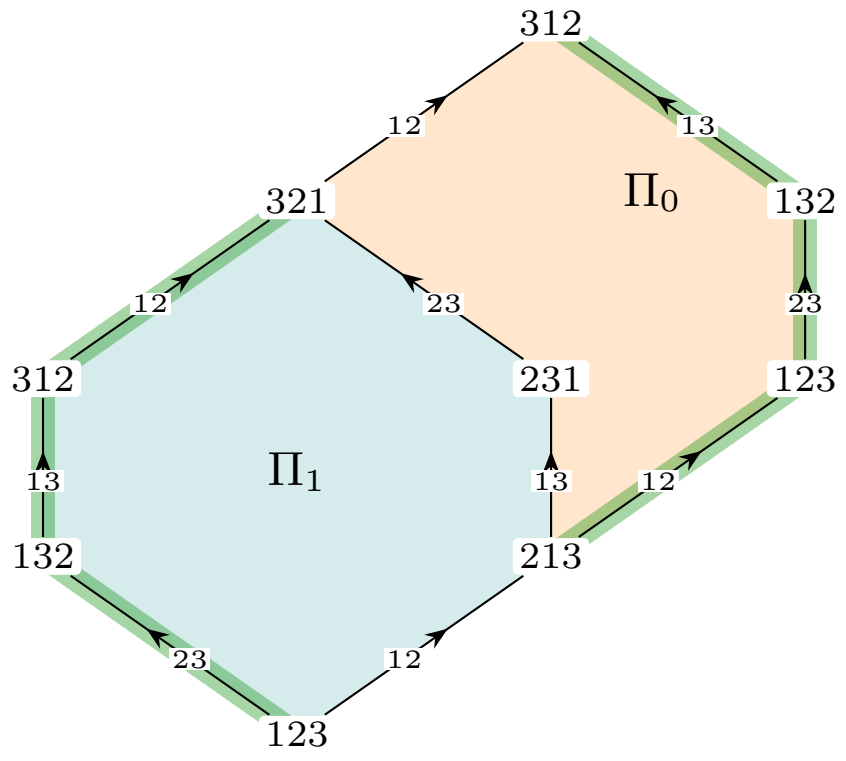
\begin{figure}[ht]
\centering
\begin{tikzpicture}[scale=1.12,vtx/.style={font=\scriptsize,inner sep=1.2pt,fill=white,rounded corners=1pt},e/.style={postaction={decorate},decoration={markings,mark=at position 0.62 with {\arrow{Stealth[length=1.5mm]}}}},lb/.style={pos=0.4,fill=white,inner sep=0.4pt,font=\tiny}]
\fill[teal!16] (1.50,0.00) -- (3.00,1.05) -- (3.00,2.10) -- (1.50,3.15) -- (0.00,2.10) -- (0.00,1.05) -- cycle;
\fill[orange!20] (3.00,1.05) -- (4.50,2.10) -- (4.50,3.15) -- (3.00,4.20) -- (1.50,3.15) -- (3.00,2.10) -- cycle;
\draw[line width=4.5pt,green!55!black,opacity=0.35,line cap=round,line join=round] (1.50,0.00) -- (0.00,1.05) -- (0.00,2.10) -- (1.50,3.15);\draw[line width=4.5pt,green!55!black,opacity=0.35,line cap=round,line join=round] (3.00,1.05) -- (4.50,2.10) -- (4.50,3.15) -- (3.00,4.20);
\node[vtx] (g) at (1.50,0.00) {123};
\node[vtx] (gb) at (0.00,1.05) {132};
\node[vtx] (gba) at (0.00,2.10) {312};
\node[vtx] (y) at (1.50,3.15) {321};
\node[vtx] (ga) at (3.00,1.05) {213};
\node[vtx] (z) at (3.00,2.10) {231};
\node[vtx] (t) at (3.00,4.20) {312};
\node[vtx] (gaa) at (4.50,2.10) {123};
\node[vtx] (gaab) at (4.50,3.15) {132};
\draw[e] (g) -- (gb) node[lb] {23};
\draw[e] (gb) -- (gba) node[lb] {13};
\draw[e] (gba) -- (y) node[lb] {12};
\draw[e] (g) -- (ga) node[lb] {12};
\draw[e] (ga) -- (z) node[lb] {13};
\draw[e] (z) -- (y) node[lb] {23};
\draw[e] (y) -- (t) node[lb] {12};
\draw[e] (ga) -- (gaa) node[lb] {12};
\draw[e] (gaa) -- (gaab) node[lb] {23};
\draw[e] (gaab) -- (t) node[lb] {13};
\node[font=\small] at (1.5,1.55) {$\Pi_1$};\node[font=\small] at (3.6,3.2) {$\Pi_0$};\node[font=\footnotesize,above=1pt] at (t.north) {peak (tip)};\node[font=\footnotesize,below=1pt] at (g.south) {valley (tip)};
\end{tikzpicture}
\caption{The diagram of Example~\ref{ex:twofaces}, with each vertex drawn at its height.
Both segments of the boundary are full: they run from the valley $123$ at the bottom up
to the peak $312$ at the top, and both ends are tips. The shaded paths are exposed
half-faces, the left side of $\Pi_1$ and the right side of $\Pi_0$.}
\label{fig:segments}
\end{figure}

\begin{example}\label{ex:twofaces}
Figure~\ref{fig:segments} shows a diagram for the word
\[
w=\sigma_2\sigma_1\sigma_2\,\sigma_2\,\sigma_1^{-1}\sigma_2^{-1}\sigma_1^{-1}\,\sigma_1^{-1},
\]
read clockwise from the valley $123$. It has two faces, $\Pi_0$ and $\Pi_1$, which share
two edges; their common vertex $231$ is an interior vertex of degree $2$, with $c=2$.
Since $w$ is longer than $R$, its area is $2$, so the diagram has the least number of
faces and is reduced. The boundary has a single peak, $312$, and a single valley, $123$.
Both are tips, and every slope has $s(x)=1$, so $\sum_{x\in\partial B}(1-s(x))=2$ and
both segments are full.

Follow the argument of Lemma~\ref{lem:exposed} down the left segment, starting from the
peak, which is the sink of $\Pi_0$. The next vertex, $321$, has degree $3$: it carries a
side corner of $\Pi_0$ and the sink corner of $\Pi_1$, so the argument passes to
$\Pi_1$. Below $321$ the whole left side of $\Pi_1$ lies on the boundary, with inner
vertices of degree $2$: this is an exposed half-face, labelled $\sigma_2\sigma_1\sigma_2$.
On the right segment the argument finds at once the exposed right side of $\Pi_0$.
Replacing $\sigma_2\sigma_1\sigma_2$ by $\sigma_1\sigma_2\sigma_1$ removes $\Pi_1$; one
deletion, a second braid move, which removes $\Pi_0$, and three more deletions then lead
to the empty word.
\end{example}

\begin{proof}[Proof of Theorem~\ref{thm:three}]
We argue by induction on the pair $(A(w),|w|)$, ordered lexicographically. If $w$ is
empty there is nothing to prove.

If some cyclic permutation of $w$ contains a subword
$\sigma_i^{\varepsilon}\sigma_i^{-\varepsilon}$, we apply this cyclic permutation and
delete the subword. The resulting word $w'$ is shorter, and $A(w')\le A(w)$ because $w'$
is freely equal to a conjugate of $w$; so $w'$ can be carried to the empty word by
induction.

Otherwise $w$ is cyclically reduced and nonempty, so $A(w)\ge1$. Let $D$ be a diagram for
$w$ with $A(w)$ faces; it is reduced, so Lemmas~\ref{lem:faces} and~\ref{lem:corners}
apply to it. No vertex of $D$ has degree $1$, for the boundary cycle of $D$ would
otherwise run along an edge and straight back, and $w$ would not be cyclically reduced.
As is standard \cite[Ch.~V]{LyndonSchupp}, $D$ is the union of its \emph{disc
components}, the closures of the components of its interior, which are disc diagrams, and
of the edges lying on no face, and these pieces are arranged like a tree: the graph with a
node for each piece and a node for each vertex lying on two or more pieces, a piece being
joined to the vertices it contains, is a tree. Since $D$ has a face, it has a disc
component. If there is more than one piece, the tree has at least two leaves, each of them
a piece meeting the rest of $D$ in a single vertex; a leaf cannot be an edge, as its other
endpoint would have degree $1$. In all cases there is therefore a disc component $B$
meeting the rest of $D$ in at most one vertex $c$.

By Lemma~\ref{lem:exposed}, $\partial B$ has two distinct segments containing exposed
half-faces whose inner vertices are slopes of the respective segments. As $c$ is a slope
of at most one segment, one of these exposed half-faces, say a path $p$ in the boundary of
a face $\Pi$, has inner vertices different from $c$. These inner vertices have degree $2$
in $B$, and hence in $D$, so the boundary cycle of $D$ traverses $p$, and some cyclic
permutation of $w$ has the form $u\,\lambda(p)$, where $\lambda(p)$ is the label of $p$.
Let $q$ be the path formed by the other three edges of $\partial\Pi$, from the initial
vertex of $p$ to its final vertex. Then $\lambda(p)\lambda(q)^{-1}$ is a cyclic
permutation of $R^{\pm1}$, so replacing $\lambda(p)$ by $\lambda(q)$ replaces one half of a
defining relator by the other. This is a braid move with non-alternating exponents: by
Lemma~\ref{lem:faces}, the heights along $p$ change by three consecutive terms of the
cyclic sequence $+1,+1,+1,-1,-1,-1$, never by $+1,-1,+1$ or $-1,+1,-1$.

Removing from $D$ the face $\Pi$, the three edges of $p$ and the two inner vertices of $p$
leaves a diagram for $w'=u\,\lambda(q)$ with $A(w)-1$ faces. Hence $A(w')<A(w)$ and
$|w'|=|w|$, and $w'$ can be carried to the empty word by induction.
\end{proof}

\begin{remark}
The proof gives more than Theorem~\ref{thm:three}: in $B_3$ every braid move in the
reduction can be chosen to lower the area by one, so a trivial word $w$ is carried to the
empty word by exactly $A(w)$ braid moves and $|w|/2$ deletions, interspersed with cyclic
permutations. The reversed inverse of a word is never needed. The contrast between three
and four strands parallels the situation for closed braids, where Morton's example
\cite{Morton1983} also needs four strands.
\end{remark}

\section{Untangling with power braid moves}\label{sec:power}

This section establishes our third contribution: the hard braids of
Section~\ref{sec:family} are not very hard, because power braid moves untangle them
without ever increasing their length.

\begin{theorem}\label{thm:power}
Under the hypotheses of Theorem~\ref{thm:hard}, and indeed for all nonzero $p,q$, the
word $w_\varepsilon(p,q)$ can be carried to the empty word by elementary moves and power
braid moves, without ever increasing its length. Exactly $2|p|$ power braid moves are
used when $|q|\ge2$.
\end{theorem}

\begin{proof}
Suppose $p>0$; the case $p<0$ is the same with $b^{-1}$ in place of $b$. In
$w_\varepsilon(p,q)=b^{-p}\,Y\,b^{p}\,Y^{-1}$, Lemma~\ref{lem:slide} slides the first
letter of $b^{p}$ to the left across $Y$ without changing the length, producing the
word $b^{-p}\,b\,Y\,b^{p-1}\,Y^{-1}$, in which $b^{-1}b$ is then deleted. After $p$ such
rounds we reach $Y\,Y^{-1}$, which is carried to the empty word by successive
deletions at its centre. Each slide uses two power braid moves, so $2p$ are used in
all.
\end{proof}

\begin{observation}
For the hard braid $W=w_1(2,2)$, Theorem~\ref{thm:power} predicts a monotone
untangling with exactly four power braid moves; an independent computer search found
a shortest such untangling, with $16$ moves of which four are power braid moves. The
power braid moves used have exponent $m=\pm2$, so they are of the even kind; in terms of the diagrams of
Section~\ref{sec:diagram}, an even power braid move is realised by two hexagons
whose tracks cross twice.
\end{observation}

\paragraph{Exponents of power braid moves.}
The proof of Theorem~\ref{thm:power} uses power braid moves of exponent $\pm q$ only,
because the slide lemma does. A power braid move of large exponent rearranges many
crossings at once, so it is natural to ask whether a fixed bound on the exponent
suffices. It does not. For $k\ge1$, let $W_k$ be the word obtained from $W$ by
replacing every $\sigma_2^{\pm2}$ by $\sigma_2^{\pm2k}$, so that $W_1=W$:
\[
  W_k=\sigma_2^{-2k}\,\sigma_3^{-1}\sigma_1\,\sigma_2^{2k}\,\sigma_1^{-1}\sigma_3\,
      \sigma_2^{2k}\,\sigma_3^{-1}\sigma_1\,\sigma_2^{-2k}\,\sigma_1^{-1}\sigma_3 .
\]
Then $W_k=w_1(2k,2k)$ with $j=2$, so $W_k$ is a trivial hard braid of length $8k+8$ by
Theorem~\ref{thm:hard}.

\begin{proposition}\label{prop:exponent}
Let $M\ge1$, and let $p,q$ be integers with $|p|>M$ and $|q|>M$. Then
$w_\varepsilon(p,q)$ cannot be carried to the empty word by elementary moves and power
braid moves of exponent $|m|\le M$, even when cyclic permutations and the passage to
the reversed inverse are allowed. In particular, $W_k$ cannot be untangled with power
braid moves of exponent less than $2k$, while by Theorem~\ref{thm:power} it can be
untangled monotonically with power braid moves of exponent $\pm2k$.
\end{proposition}

\begin{proof}
Since $|p|,|q|\ge2$, Theorem~\ref{thm:hard} applies. Let $\Omega$ be the set of words
introduced in its proof; it is closed under elementary moves, cyclic permutations and
reversed inverses, and none of its words admits a deletion. It suffices to show that
no power braid move of exponent $m$ with $2\le|m|\le M$ applies to a word of $\Omega$;
power braid moves with $|m|=1$ are elementary braid moves.

Such a move needs a subword $\sigma_i^{\varepsilon}\sigma_k^{m}\sigma_i^{-\varepsilon}$
with $|i-k|=1$, that is, $|m|\ge2$ equal consecutive letters flanked by two letters of an
adjacent generator. In a word of $\Omega$, two consecutive letters are equal only inside
a run of $b$: the two letters of a block belong to different generators, and a block
meets a run in a letter of $a$ or $c$ and a letter of $b$. Hence $\sigma_k=b$, the
flanking letters belong to $a$ or $c$, and the $|m|$ letters $b^{\pm1}$ form a whole
maximal run. Its length is $|p|$ or $|q|$, both greater than $M$, which is
impossible.
\end{proof}

So the third contribution comes with a caveat: power braid moves untangle every member
of the family, but the exponent they require grows along the family, and no bound on
the exponent of power braid moves suffices to untangle all the braids $W_k$. Whether some
other finite set of length-preserving moves untangles all of them is a question we leave
open.

\begin{remark}\label{rem:threestrands}
Power braid moves are local in the strands, though not in the crossings. The generators
$\sigma_i$ and $\sigma_k$ with $|i-k|=1$ act only on the positions $\min(i,k)$,
$\min(i,k)+1$ and $\min(i,k)+2$. So throughout either side of a power braid move, the
three strands in these positions are permuted among themselves, and the other $n-3$
strands pass straight through. Both sides lie in the subgroup generated by $\sigma_i$ and
$\sigma_k$, which is isomorphic to $B_3$, and the move is the $m$-th power of a braid
relation. In this respect a power braid move is no less local than an elementary braid
move. What it does not respect is the number of crossings: it replaces $|m|+2$
crossings at once.

This interacts with Theorem~\ref{thm:three} as follows. The word
$\sigma_i^{\varepsilon}\sigma_k^{m}\sigma_i^{-\varepsilon}\,
\sigma_k^{-\varepsilon}\sigma_i^{-m}\sigma_k^{\varepsilon}$, the left-hand side of the move
followed by the inverse of the right-hand side, is a trivial braid on three strands; so
by Theorem~\ref{thm:three} it can be untangled without ever increasing its length. This
does not mean that one side of the move can be turned into the other by elementary moves
without lengthening. For instance, no elementary move applies to the word
$\sigma_1\sigma_2^{2}\sigma_1^{-1}$ at all: it has no cancellation, no commuting pair of
adjacent letters, and no subword of the form $x\,y\,x'$ with $x,x'$ letters of the same
generator.

Thus the hard braids of Section~\ref{sec:family} are untangled by moves each of which
involves only three strands, and Proposition~\ref{prop:exponent} shows that the number
of crossings these moves must handle at once is unbounded: $W_k$ needs a power braid
move whose two sides each have $2k+2$ crossings, all among the same three strands. The elementary moves fail on $W_k$ not because they involve too few strands, but because
they involve too few crossings. This observation suggests that moves on three strands, of
unbounded length, might untangle every trivial braid; Section~\ref{sec:bands} shows that
they do not.
\end{remark}

\section{Moves on a bounded number of strands}\label{sec:bands}

This section establishes our fourth contribution. The hard braids of
Section~\ref{sec:family} resist the elementary moves, but they yield to power braid moves,
which are longer but still involve only three strands. One may hope that this is the general
picture: that every trivial braid can be untangled monotonically by moves which are allowed
to be long but are confined to a small number of strands. We show that this hope fails for
every bound on the number of strands.

\subsection{Moves on \texorpdfstring{$k$}{k} strands, and a question in the spirit of Tits}

\begin{definition}\label{def:kmove}
Let $2\le k\le n$. A subword $s$ of a braid word on $n$ strands is a \emph{$k$-strand subword}
if, for some $a$, all its letters are among
$\sigma_a^{\pm1},\dots,\sigma_{a+k-2}^{\pm1}$. A \emph{move on $k$ strands} replaces a
$k$-strand subword $s$ by a word $s'$ in the same letters
$\sigma_a^{\pm1},\dots,\sigma_{a+k-2}^{\pm1}$ such that
\[
  s'=s\ \text{ in } B_n \qquad\text{and}\qquad |s'|\le|s| .
\]
The move is \emph{inessential} if $s'$ can be obtained from $s$ by commutations (elementary
moves of type~2), and \emph{essential} otherwise. A trivial braid word can be
\emph{untangled monotonically by moves on at most $k$ strands} if it can be carried to the
empty word by moves on at most $k$ strands, commutations, cyclic permutations and passages
to the reversed inverse.
\end{definition}

A $k$-strand subword involves only the $k$ strands that occupy the positions
$a,\dots,a+k-1$; the other strands pass straight through. So a move on $k$ strands replaces a
portion of the braid in which only $k$ adjacent strands take part by another braid on the same
strands, equal to it and with no more crossings. There is no bound on the length of $s$, so
for $k\ge3$ there are infinitely many moves on $k$ strands. The subgroup of $B_n$ generated by
$\sigma_a,\dots,\sigma_{a+k-2}$ is a copy of $B_k$, so whether a given replacement is a move
can be decided inside $B_k$. A move on $k$ strands is also a move on $k'$ strands whenever
$k\le k'\le n$. Commutations are allowed freely in all that follows: a commutation exchanges
two crossings which have no strand in common, and does not change the braid diagram up to
planar isotopy.

Deletions are moves on $2$ strands. Elementary braid moves and power braid moves are moves on
$3$ strands (Remark~\ref{rem:threestrands}). Moves on $3$ strands are closely related to
the following moves, suggested by Tits' theorem.

\begin{definition}
A \emph{two-generator move} replaces a subword $s$ of a braid word, in which only two
generators $\sigma_i,\sigma_k$ occur, by a word $t$ in $\sigma_i^{\pm1},\sigma_k^{\pm1}$
such that $t=s$ in $B_n$ and $|t|\le|s|$.
\end{definition}

If $|i-k|=1$, a two-generator move is a move on $3$ strands, and every move on $3$ strands
arises in this way. If $|i-k|\ge2$, the subgroup generated by $\sigma_i$ and $\sigma_k$ is
free abelian of rank $2$, and a two-generator move is a composite of commutations and
deletions. Hence a trivial braid word can be untangled monotonically by two-generator moves
if and only if it can be untangled monotonically by moves on at most $3$ strands.

Tits' theorem, recalled in the introduction, can be read as a statement about such moves. In
a Coxeter group, every trivial word can be carried to the empty word, without ever increasing
its length, by deletions $ss$ and braid relations, and each of these moves takes place inside
a subgroup generated by at most two of the standard generators. The hard braids of
Section~\ref{sec:family} show that for braid groups the braid relations themselves are not
enough; but all the moves used in Theorem~\ref{thm:power} still take place inside subgroups
generated by two standard generators. This suggests the following analogue of Tits' theorem,
in which the braid relations are replaced by \emph{all} the length-non-increasing relations of
these subgroups. For $k\ge2$, consider the statement
\begin{itemize}
\item[$(\mathrm T_k)$] \emph{every trivial braid word, on any number of strands, can be
  untangled monotonically by moves on at most $k$ strands.}
\end{itemize}
Statement $(\mathrm T_3)$ says that two-generator moves always suffice; an earlier version
of this article proposed it as a conjecture. There was reasonable evidence for it:
\begin{itemize}
\item On three strands it holds trivially, since the whole word is a two-generator word; and
  Theorem~\ref{thm:three} shows that even elementary moves suffice.
\item If a trivial word has the form $v_1v_2^{-1}$ with $v_1$ and $v_2$ positive words, then
  $v_1=v_2$ in the positive braid monoid by Garside's theorem \cite{Garside}, so $v_1$ can be
  transformed into $v_2$ by positive braid relations, which are elementary moves, after which
  $v_2v_2^{-1}$ is deleted letter by letter.
\item The family of Section~\ref{sec:family}, including the braids $W_k$, satisfies it by
  Theorem~\ref{thm:power}, and so does every hard braid that our searches found
  (Section~\ref{sec:verification}).
\item An exhaustive search shows that $(\mathrm T_3)$ holds for all short words: every
  trivial braid word of length at most $16$ in $B_4$, at most $14$ in $B_5$ and at most $12$
  in $B_6$ can be untangled monotonically by moves on at most $3$ strands
  (Section~\ref{sec:verification}).
\item For Artin groups of large type, in which every pair of standard generators satisfies a
  braid relation of length at least $3$, the analogue of $(\mathrm T_3)$ is a theorem of Holt
  and Rees \cite[Theorem~3.2]{HoltRees}. They show that every word can be carried to its
  shortlex normal form, without ever increasing its length, by free cancellations and
  $\tau$-moves, which replace a two-generator subword by another word of the same length
  representing the same element; in particular every trivial word is carried to the empty
  word. Their argument uses the largeness assumption, and braid groups on four or more
  strands are not of large type, because $\sigma_i$ and $\sigma_k$ commute when $|i-k|\ge2$.
\end{itemize}
Nevertheless $(\mathrm T_3)$ is false, and so is $(\mathrm T_k)$ for every $k$: this is
Theorem~\ref{thm:bands} below. Before turning to it, we record that, like power braid moves
of bounded exponent, two-generator moves of bounded length do not even untangle the family of
Section~\ref{sec:family}.

\begin{proposition}\label{prop:boundedlength}
Let $M\ge3$, and let $p,q$ be integers with $|p|\ge M-1$ and $|q|\ge M-1$. Then
$w_\varepsilon(p,q)$ cannot be carried to the empty word by two-generator moves on subwords
of length at most $M$, even when cyclic permutations and the passage to the reversed
inverse are allowed. In particular, $W_k$ requires two-generator moves on subwords of length
$2k+2$, and by Theorem~\ref{thm:power} such moves suffice.
\end{proposition}

\begin{proof}
Since $|p|,|q|\ge2$, Theorem~\ref{thm:hard} applies; let $\Omega$ be the set of words
introduced in its proof. It is closed under elementary moves, cyclic permutations and
reversed inverses, and none of its words admits a deletion. We show that every
two-generator move on a subword $s$ of length at most $M$ of a word in $\Omega$ either does
nothing or is a commutation inside a block, so that $\Omega$ is closed under these moves too;
the empty word is then unreachable.

\emph{The possible subwords.} In a word of $\Omega$, the letters of $a$ and $c$ occur only
in the blocks, each block consists of one letter of $a$ and one of $c$, and consecutive
blocks are separated by a maximal run of letters $b^{\pm1}$ of length $|p|$ or $|q|$. Hence
a subword $s$ in which only two generators occur either lies inside a block, and then any
two-generator move on it is a commutation or does nothing, or it involves only $b$ and one of
$a$, $c$. In the latter case $s$ consists of part of a run $b^{\pm r}$, possibly preceded by
the letter $x$ of the preceding block adjacent to the run and possibly followed by the letter
$y$ of the following block adjacent to the run. If both $x$ and $y$ occur, $s$ contains the
whole run, so that $|s|\ge\min(|p|,|q|)+2>M$. So $s$ has the form $b^{\pm r'}$,
$x\,b^{\pm r'}$ or $b^{\pm r'}y$.

\emph{These subwords are unique geodesics.} We use the description of geodesics in
$B_3=\langle x,y\mid xyx=yxy\rangle$ due to Mairesse and Math\'eus \cite{MairesseMatheus},
as stated in \cite[Section~2]{HoltRees}. For a freely reduced word $w$, let $p(w)$ and $n(w)$
be the lengths of its longest subwords consisting of alternating positive letters and of
alternating negative letters, respectively, capped at $3$. Then $w$ is geodesic if and only
if $p(w)+n(w)\le3$, and it is the unique geodesic word representing its element if
$p(w)+n(w)<3$. For $s=b^{\pm r'}$ we have $p(s)+n(s)=1$. For $s=x\,b^{\pm r'}$ or
$b^{\pm r'}y$ we have $p(s)+n(s)=2$: the longest alternating subwords have length $2$ and
$0$ if the two letters involved have the same sign, and length $1$ each otherwise. In all
three cases $s$ is the unique geodesic representing its element in the subgroup generated by
its two generators, which is isomorphic to $B_3$; so any word $t$ with $t=s$ and $|t|\le|s|$
equals $s$, and the move does nothing.

For $W_k=w_1(2k,2k)$ we may take $M=2k+1$, and the untangling of Theorem~\ref{thm:power} uses
power braid moves, which are two-generator moves on subwords of length $2k+2$.
\end{proof}

\subsection{Band braids}\label{subsec:bandbraids}

The braids that defeat moves on few strands are obtained from the hard braid $w_{-1}(-2,2)$
of Section~\ref{sec:family} by replacing each of its four strands by a \emph{band} of $m$
parallel strands.

For integers $a,b\ge0$ and $q\ge0$, put
\[
  \chi_q(a,b)=\prod_{i=a}^{1}\ \prod_{j=1}^{b}\sigma_{q+i+j-1}
  =(\sigma_{q+a}\sigma_{q+a+1}\cdots\sigma_{q+a+b-1})\,
   (\sigma_{q+a-1}\cdots\sigma_{q+a+b-2})\cdots(\sigma_{q+1}\cdots\sigma_{q+b}),
\]
the outer product being taken over $i=a,a-1,\dots,1$ in this order, and let
$\bar\chi_q(a,b)$ be the word obtained from $\chi_q(a,b)$ by inverting every letter, without
changing their order. In $\chi_q(a,b)$ a band of $a$ strands, in positions
$q+1,\dots,q+a$, moves to the right \emph{over} a band of $b$ strands, in positions
$q+a+1,\dots,q+a+b$; the strands of each band stay parallel. The word has $ab$ letters, one
for each pair consisting of a strand of the left band and a strand of the right band. Reading
the definition, the $i$-th strand of the left band meets the strands of the right band in the
order $1,2,\dots,b$, and the $j$-th strand of the right band meets the strands of the left
band in the order $a,a-1,\dots,1$. In $\bar\chi_q(a,b)$ the left band moves to the right
\emph{under} the right band, with the same orders. For instance
$\chi_2(2,2)=\sigma_4\sigma_5\sigma_3\sigma_4$.

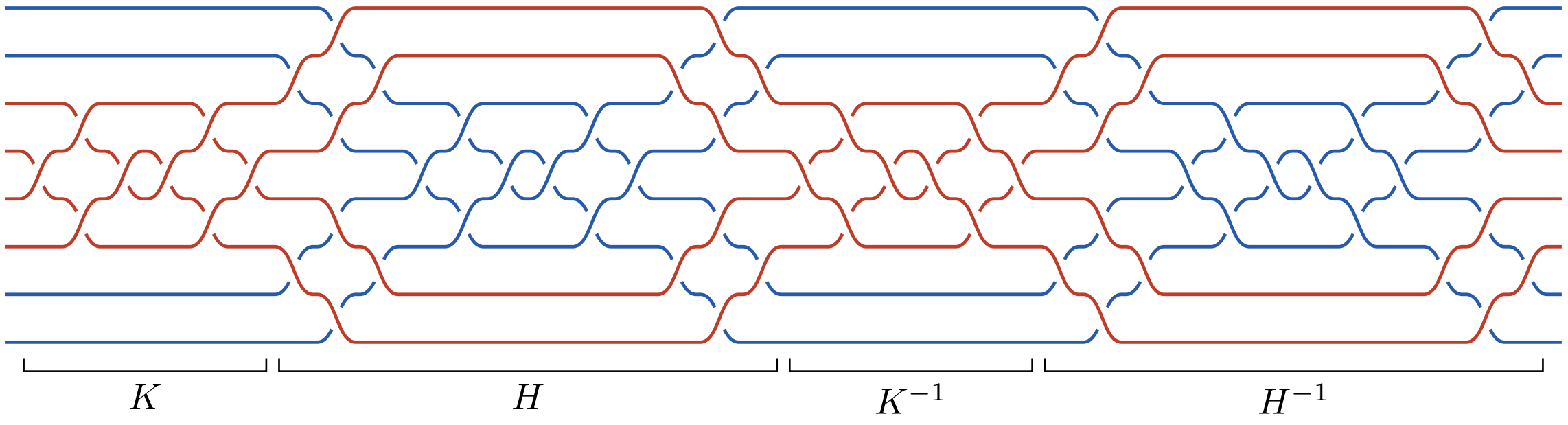
\begin{figure}[ht]
\centering
\rotatebox{90}{%
\begin{tikzpicture}
\pic[braid/.cd, number of strands=8, width=4.6mm, height=-4.1mm, border height=1.2mm, gap=0.16,
     every strand/.style={line width=0.9pt},
     strand 1/.style={botc}, strand 2/.style={botc}, strand 7/.style={botc}, strand 8/.style={botc},
     strand 3/.style={topc}, strand 4/.style={topc}, strand 5/.style={topc}, strand 6/.style={topc}]
 {braid={ s_4 s_3-s_5 s_4 s_4 s_3-s_5 s_4
          s_2^{-1}-s_6 s_1^{-1}-s_3^{-1}-s_5-s_7 s_2^{-1}-s_6
          s_4 s_3-s_5 s_4 s_4 s_3-s_5 s_4
          s_6^{-1}-s_2 s_7^{-1}-s_5^{-1}-s_1-s_3 s_6^{-1}-s_2
          s_4^{-1} s_5^{-1}-s_3^{-1} s_4^{-1} s_4^{-1} s_5^{-1}-s_3^{-1} s_4^{-1}
          s_2^{-1}-s_6 s_1^{-1}-s_3^{-1}-s_5-s_7 s_2^{-1}-s_6
          s_4^{-1} s_5^{-1}-s_3^{-1} s_4^{-1} s_4^{-1} s_5^{-1}-s_3^{-1} s_4^{-1}
          s_6^{-1}-s_2 s_7^{-1}-s_5^{-1}-s_1-s_3 s_6^{-1}-s_2 }};
\end{tikzpicture}}\\[1mm]
\begin{tikzpicture}[x=4.1mm]
\foreach \a/\b/\l in {0/6/$K$, 6/18/$H$, 18/24/$K^{-1}$, 24/36/$H^{-1}$}{
  \draw[line width=0.5pt] (\a+0.15,0.12) -- (\a+0.15,0) -- (\b-0.15,0) -- (\b-0.15,0.12);
  \node[below] at ({(\a+\b)/2},0) {\l};}
\end{tikzpicture}
\caption{The band braid $V_2=K\,H\,K^{-1}H^{-1}$ on $8$ strands, with $64$ crossings, read
from left to right, with position~$1$ at the bottom. The upper bands $\mathcal P$ and
$\mathcal Q$ are drawn in red and the lower bands $\mathcal P'$ and $\mathcal Q'$ in blue; a
red strand always passes over a blue one. Crossings drawn in the same column commute. No
move on at most four strands changes this word, except by commutations
(Theorem~\ref{thm:bands}).}
\label{fig:band}
\end{figure}

Now fix $m\ge1$ and work with $4m$ strands, divided into four bands of $m$ strands each,
\[
  \mathcal P'=\{P'_1,\dots,P'_m\},\quad \mathcal P=\{P_1,\dots,P_m\},\quad
  \mathcal Q=\{Q_1,\dots,Q_m\},\quad \mathcal Q'=\{Q'_1,\dots,Q'_m\},
\]
which initially occupy the positions $1,\dots,m$; $m+1,\dots,2m$; $2m+1,\dots,3m$ and
$3m+1,\dots,4m$, in this order. We call $\mathcal P$ and $\mathcal Q$ the \emph{upper} bands
and $\mathcal P'$ and $\mathcal Q'$ the \emph{lower} bands. Put
\[
  K=\chi_m(m,m)^2,\qquad E=\bar\chi_0(m,m),\qquad F=\chi_{2m}(m,m),\qquad
  H=E\,F\,K\,F^{-1}E^{-1},
\]
and
\[
  V_m=K\,H\,K^{-1}H^{-1}.
\]
Here is what these words do (Figure~\ref{fig:band}). The word $\chi_m(m,m)$ makes the band
in positions $m+1,\dots,2m$ pass over the band in positions $2m+1,\dots,3m$, and $K$ does
this twice, so that the two bands are back in place: we say that the two bands form a
\emph{clasp}. In $V_m$ the first factor $K$ is a clasp of the upper bands $\mathcal P$ and
$\mathcal Q$. In $H$, the word $E$ slides $\mathcal P'$ to the right underneath $\mathcal P$,
and $F$ slides $\mathcal Q'$ to the left underneath $\mathcal Q$; the bands are then in the
order $\mathcal P,\mathcal P',\mathcal Q',\mathcal Q$, so that the factor $K$ of $H$ is a
clasp of the lower bands $\mathcal P'$ and $\mathcal Q'$; finally $F^{-1}$ and $E^{-1}$ bring
the lower bands back. The words $K$ and $H$ have $2m^2$ and $6m^2$ letters, so $V_m$ has
length $16m^2$. In every crossing between an upper strand and a lower strand, the upper strand
passes over the lower one, and the lower bands clasp \emph{between} the two upper bands, which
are themselves clasped. To cancel the two clasps against their inverses, one of them has to be
moved out of the way of the other, and Theorem~\ref{thm:bands} says how many strands this
takes.

For $m=1$ the bands are single strands and
\[
  V_1=\sigma_2^{2}\cdot\sigma_1^{-1}\sigma_3\,\sigma_2^{2}\,\sigma_3^{-1}\sigma_1\cdot
      \sigma_2^{-2}\cdot\sigma_1^{-1}\sigma_3\,\sigma_2^{-2}\,\sigma_3^{-1}\sigma_1
\]
is the hard braid $w_{-1}(-2,2)$ of Theorem~\ref{thm:hard}, with $j=2$. For $m=2$, writing $i$
for $\sigma_i$ and $\bar\imath$ for $\sigma_i^{-1}$,
\begin{align*}
K&=4\,5\,3\,4\;4\,5\,3\,4,\\
H&=\bar2\,\bar3\,\bar1\,\bar2\;\,6\,7\,5\,6\;\,4\,5\,3\,4\;4\,5\,3\,4\;\,
   \bar6\,\bar5\,\bar7\,\bar6\;\,2\,1\,3\,2,
\end{align*}
and $V_2=K\,H\,K^{-1}H^{-1}$ is a braid word on $8$ strands with $64$ crossings.

\begin{theorem}\label{thm:bands}
Let $m\ge1$ and $n\ge4m$, and regard $V_m$ as a braid word on $n$ strands.
\begin{enumerate}
\item[\textup{(i)}] $V_m$ is a trivial braid word of length $16m^2$.
\item[\textup{(ii)}] Let $w$ be any word obtained from $V_m$ by commutations, cyclic
  permutations and passages to the reversed inverse. Then every move on at most $2m$ strands
  that can be applied to $w$ is inessential: it replaces a subword $s$ by a word that
  can be obtained from $s$ by commutations. Consequently $V_m$ cannot be untangled monotonically
  by moves on at most $2m$ strands; indeed, such moves never shorten it.
\item[\textup{(iii)}] $V_m$ can be untangled monotonically by moves on $2m+1$ strands: $4m$
  such moves, all of them length-preserving, together with commutations and deletions, carry
  $V_m$ to the empty word.
\end{enumerate}
\end{theorem}

For $m=1$, part~(ii) says only that $V_1$ admits no deletion, and part~(iii) recovers the
untangling of $V_1$ by four power braid moves given by Theorem~\ref{thm:power}. The new
content starts at $m=2$.

\begin{corollary}\label{cor:tits}
\begin{enumerate}
\item[\textup{(a)}] The braid word $V_2$, on $8$ strands and with $64$ crossings, cannot be
  untangled monotonically by two-generator moves. Hence statement $(\mathrm T_3)$ is false
  in $B_n$ for every $n\ge8$. In particular, elementary moves and power braid moves, of
  arbitrary exponents, do not untangle $V_2$.
\item[\textup{(b)}] Statement $(\mathrm T_k)$ is false for every $k$: the trivial braid word
  $V_{\lceil k/2\rceil}$ cannot be untangled monotonically by moves on at most $k$ strands.
\item[\textup{(c)}] For $n\ge2$ let $\kappa(n)$ be the least $k$ such that every trivial
  braid word on $n$ strands can be untangled monotonically by moves on at most $k$ strands.
  Then
  \[
     2\lfloor n/4\rfloor+1\ \le\ \kappa(n)\ \le\ n .
  \]
\end{enumerate}
\end{corollary}

\begin{proof}
(a) and (b) follow from Theorem~\ref{thm:bands}(ii), with $m=2$ and with
$m=\lceil k/2\rceil$, since a two-generator move is either a move on $3$ strands or a
composite of commutations and deletions, and elementary and power braid moves are
two-generator moves. For
(c), the upper bound holds because the replacement of a whole trivial word by the empty word
is a move on $n$ strands, and the lower bound is Theorem~\ref{thm:bands}(ii) with
$m=\lfloor n/4\rfloor$ when $n\ge4$, and is empty when $n<4$.
\end{proof}

Thus braid groups on many strands behave, in this respect, like diagrams of knots and not
like Coxeter groups or Artin groups of large type: no family of length-non-increasing moves
confined to a bounded number of strands untangles all trivial braids. The proof of
Theorem~\ref{thm:bands} occupies the rest of this section. Parts~(i) and~(iii) are proved in
Section~\ref{subsec:layers} by a simple geometric principle. Part~(ii) is proved in
Sections~\ref{subsec:itineraries}--\ref{subsec:proofii}; for $m=2$ and $m=3$ it was also
checked by exhaustive computer searches, described in Section~\ref{sec:verification}, which
do not rely on the proof.

\paragraph{Outline of the proof of part~(ii).}
Suppose that a subword $s$ involving at most $2m$ adjacent strands is replaced by an equal
braid word $s'$ with no more crossings. We must show that $s'$ is $s$ up to commutations. The
proof has two steps.

\emph{Step 1: $s'$ has the same number of crossings as $s$ between any two strands.} A braid
word is as long as the sum, over all pairs of strands, of the number of times the two strands
cross. We bound these numbers from below in $s'$ by looking at two or three strands at a time.
The key observation is that $V_m$ restricts to a copy of the hard braid $V_1$ on any four
strands taken one from each band, and that, because at most $2m$ strands are involved, there
are always enough strands \emph{outside} $s$ to choose such copies of $V_1$ which meet $s$ in
only two or three strands. On two or three strands the relevant parts of $V_1$ are as short
as possible, and a short counting argument assembles these local bounds into the equality
of all the crossing numbers.

\emph{Step 2: the crossings occur in the same order.} Once the numbers of crossings are
known, the order in which each strand meets the others is forced: two strands of the same
band never cross, so a strand can only reach a strand of another band by crossing those in
between. A braid word is determined, up to commutations, by the sequence of crossings along
each strand (Lemma~\ref{lem:itin}), and this finishes the proof.

Sections~\ref{subsec:itineraries}--\ref{subsec:three} prepare the tools: the description
of a braid word by the itineraries of its strands, the itineraries of $V_m$, and two lemmas
on which configurations can occur among at most $2m$ strands.

\subsection{Layered braids: proof of parts (i) and (iii)}\label{subsec:layers}

Throughout this section we use the language of braid diagrams. We name the strands of a braid
word by their initial positions. Each letter $\sigma_i^{\pm1}$ is a \emph{crossing} of two
strands, namely those which occupy the positions $i$ and $i+1$ when the letter is read, and
it has a \emph{sign}, its exponent; as in Section~\ref{sec:diagram}, in a positive crossing the
strand on the left passes over the strand on the right. For a set $T$ of strands of a braid word
$w$, we write $w|_T$ for the braid word on $|T|$ strands obtained by deleting all the other
strands, that is, by keeping only the crossings between two strands of $T$. If $w=w'$ in
$B_n$ then $w|_T=w'|_T$ in $B_{|T|}$, since an isotopy of a geometric braid restricts to an
isotopy of any subset of its strands.

Suppose the strands of a braid word are divided into two classes, \emph{upper} and
\emph{lower}. We say that the word is \emph{layered} if at each crossing between an upper
strand and a lower strand, the upper strand passes over the lower one.

\begin{lemma}[layers]\label{lem:layer}
Let $w$ and $w'$ be braid words on the same strands, inducing the same permutation, and
layered with respect to the same division of the strands into a set $U$ of upper strands and a set
$L$ of lower strands. If $w|_U=w'|_U$ and $w|_L=w'|_L$ as braids, then $w=w'$ as braids.
\end{lemma}

\begin{proof}
Realise the braids geometrically in $\mathbb R^2\times[0,1]$, the first coordinate of
$\mathbb R^2$ being the horizontal position in the diagram and the second the height above
the plane of the diagram; a braid is a collection of disjoint strands, each meeting every
level $\mathbb R^2\times\{t\}$ in one point, considered up to isotopy through such
collections with fixed endpoints. Moving the endpoints of the upper strands up to height
$1$ and those of the lower strands down to height $-1$, by one and the same motion at both
ends of the braid, does not affect whether two braids with the same endpoints are isotopic.
Since $w$ is layered, its diagram is then the projection of a geometric braid whose upper
strands lie entirely in the half-space of positive height and whose lower strands lie
entirely in the half-space of negative height; and likewise for $w'$. The upper strands of
$w$ and those of $w'$ form two braids in the upper half-space, with the same endpoints
because $w$ and $w'$ induce the same permutation, and they are isotopic because
$w|_U=w'|_U$; as the half-space is a copy of $\mathbb R^2\times[0,1]$, the isotopy can be
performed inside it. The same holds for the lower strands in the lower half-space. The two
isotopies have disjoint supports, so together they form an isotopy from $w$ to $w'$.
\end{proof}

\begin{proof}[Proof of Theorem~\ref{thm:bands}\textup{(i)}]
Let $U=\mathcal P\cup\mathcal Q$ and $L=\mathcal P'\cup\mathcal Q'$. The crossings of $V_m$
between an upper and a lower strand are those of the four subwords $E^{\pm1}$ and
$F^{\pm1}$. In $E=\bar\chi_0(m,m)$ the band $\mathcal P'$, on the left, passes under
$\mathcal P$, and $E^{-1}$ undoes this; in $F=\chi_{2m}(m,m)$ the band $\mathcal Q$, on the
left, passes over $\mathcal Q'$, and $F^{-1}$ undoes this. So $V_m$ is layered. Deleting the lower
strands leaves the word $K\,K^{-1}$ on the strands of $U$, and deleting the upper strands
leaves $K\,K^{-1}$ on the strands of $L$, both trivial. As $V_m$ induces the identity
permutation, Lemma~\ref{lem:layer}, applied to $V_m$ and the empty word, shows that
$V_m$ is trivial.
\end{proof}

A typical use of Lemma~\ref{lem:layer} is this: if one strand passes over all the other
strands it meets, it may be moved freely across the braid formed by the other strands,
provided its endpoints are kept. All the moves in the following proof are of this kind.

\begin{proof}[Proof of Theorem~\ref{thm:bands}\textup{(iii)}]
We use repeatedly that a band crossing can be performed one strand at a time: up to
commutations,
\[
  \chi_q(a,b)\equiv\chi_q(a,b-1)\,\chi_{q+b-1}(a,1),\qquad
  \chi_q(a,b)=\chi_{q+1}(a-1,b)\,\chi_q(1,b),
\]
and similarly for $\bar\chi$; here and below $\equiv$ means that two words differ by
commutations.

\emph{Stage 1: the upper strands of $\mathcal P$ step across the lower clasp.} For
$0\le j\le m$ put
\[
  H_j=E_j\,F\,G_j\,K_j\,G_j^{-1}F^{-1}E_j^{-1},\qquad
  E_j=\bar\chi_0(m,j),\quad G_j=\chi_{m+j}(m-j,m),\quad K_j=\chi_j(m,m)^2 .
\]
In $H_j$ the band $\mathcal P'$ passes under $P_1,\dots,P_j$ only, while $\mathcal Q'$
passes under $\mathcal Q$ and then under $P_{j+1},\dots,P_m$; the two lower bands then meet
between $P_j$ and $P_{j+1}$, where they clasp. Thus $H_m=H$, and each $H_j$ has $6m^2$ letters.
Let $1\le j\le m$, and let $e_j=\bar\chi_{j-1}(m,1)$ and $d_j=\chi_{m+j-1}(1,m)$. Then
$E_j\equiv E_{j-1}e_j$ and $G_{j-1}=G_j\,d_j$, and $e_j$ commutes with $F$ and with $G_j$,
whose letters have indices greater by at least $2$. Hence
\[
  H_j\ \equiv\ E_{j-1}\,F\,G_j\,\bigl[e_j\,K_j\,e_j^{-1}\bigr]\,G_j^{-1}F^{-1}E_{j-1}^{-1}.
\]
The bracketed subword involves the $2m+1$ strands of $\mathcal P'\cup\{P_j\}\cup\mathcal Q'$,
in positions $j,\dots,j+2m$: the strand $P_j$ steps to the left over $\mathcal P'$, the lower
bands clasp, and $P_j$ steps back. The word $d_j\,K_{j-1}\,d_j^{-1}$ involves the same
strands: now $P_j$ steps to the right over $\mathcal Q'$, the lower bands clasp, and $P_j$
steps back. Both words are layered, with $P_j$ as the only upper strand, both induce the
identity permutation, and both reduce to a clasp of $\mathcal P'$ and $\mathcal Q'$ when
$P_j$ is deleted; so they are equal by Lemma~\ref{lem:layer}. Both have $2m^2+2m$ letters.
Replacing the first by the second is therefore a length-preserving move on $2m+1$ strands,
and it turns $H_j$ into $E_{j-1}F\,G_j\,d_j\,K_{j-1}\,d_j^{-1}G_j^{-1}F^{-1}E_{j-1}^{-1}=H_{j-1}$.
After $m$ such moves $H$ has become $H_0$; the inverse moves turn $H^{-1}$ into
$H_0^{-1}$. So $2m$ moves carry $V_m$ to $K\,H_0\,K^{-1}H_0^{-1}$.

\emph{Stage 2: the lower strands of $\mathcal Q'$ slide under the upper clasp.} We have
$H_0=M\,K_0\,M^{-1}$ with $M=F\,G_0$, the word in which $\mathcal Q'$ passes to the left
under $\mathcal Q$ and then under $\mathcal P$. Up to commutations,
$M\equiv\gamma_1\gamma_2\cdots\gamma_m$ with $\gamma_l=\chi_{m+l-1}(2m,1)$, the word in which
the single strand $Q'_l$ passes under the $2m$ upper strands. Let
$K^{(l)}=\chi_{m+l}(m,m)^2$ be the clasp of the upper bands when they occupy the positions
$m+l+1,\dots,3m+l$, so that $K^{(0)}=K$. For $1\le l\le m$ the words
$K^{(l-1)}\gamma_l$ and $\gamma_l\,K^{(l)}$ involve the same $2m+1$ strands, are layered
with $Q'_l$ as the only lower strand, induce the same permutation, and reduce to the clasp of
$\mathcal P$ and $\mathcal Q$ when $Q'_l$ is deleted; so they are equal by
Lemma~\ref{lem:layer}, and they have the same length. Applying these $m$ moves in turn
carries $K\,M$ to $M\,K^{(m)}$, and in the same way $m$ moves carry $K^{-1}M$ to
$M\,(K^{(m)})^{-1}$. Hence $2m$ further moves on $2m+1$ strands carry
\[
  K\,M\,K_0\,M^{-1}\,K^{-1}\,M\,K_0^{-1}M^{-1}
  \quad\text{to}\quad
  M\,K^{(m)}\,K_0\,M^{-1}M\,(K^{(m)})^{-1}K_0^{-1}M^{-1}.
\]
Deletions remove $M^{-1}M$. The letters of $K^{(m)}$ have indices at least $2m+1$ and those
of $K_0$ have indices at most $2m-1$, so commutations and deletions remove
$K^{(m)}K_0(K^{(m)})^{-1}K_0^{-1}$, and finally $MM^{-1}$.
\end{proof}

\subsection{Itineraries and convex sets of crossings}\label{subsec:itineraries}

The proof of part~(ii) describes a braid word by what happens along each of its strands.
The \emph{itinerary} of a strand $x$ in a braid word $w$ is the sequence of the crossings of
$w$ on $x$, in the order in which they occur, each recorded with the other strand it involves
and with its sign. Commutations do not change itineraries, because two letters which are
exchanged by a commutation have no strand in common. The converse also holds.

\begin{lemma}\label{lem:itin}
Let $w$ and $w'$ be braid words on the same strands, in the same initial positions, such
that every strand has the same itinerary in $w$ and in $w'$. Then $w'$ can be obtained from $w$ by commutations.
\end{lemma}

\begin{proof}
The two words have the same length, half the total length of the itineraries; we argue by
induction on it. Let the first letter of $w$ be $\sigma_i^{\varepsilon}$, a crossing of the
strands $x$ and $y$ which initially occupy the positions $i$ and $i+1$. It is the first entry
of the itineraries of $x$ and of $y$. So in $w'$ the first crossing on $x$ is a crossing with
$y$ of sign $\varepsilon$, and it is also the first crossing on $y$. Before it, neither $x$
nor $y$ has moved, so this crossing is a letter $\sigma_i^{\varepsilon}$, and the letters
preceding it in $w'$ involve neither position $i$ nor position $i+1$: they are letters
$\sigma_j^{\pm1}$ with $|j-i|\ge2$. Commutations therefore bring this letter to the front of
$w'$. Removing the first letter from both words leaves two words with equal itineraries, and
we conclude by induction.
\end{proof}

For two crossings $\alpha,\beta$ of a braid word $w$, write $\alpha\prec\beta$ if there are
crossings $\alpha=\gamma_0,\gamma_1,\dots,\gamma_r=\beta$ with $r\ge1$, such that each
$\gamma_{i}$ occurs later in $w$ than $\gamma_{i-1}$ and has a strand in common with it. This is
a partial order, unchanged by commutations; it is the usual order on the letters of a word
in a partially commutative monoid \cite{CartierFoata}. If $\alpha$ and $\beta$ lie on a common strand, the
one occurring first precedes the other. A set $S$ of crossings is \emph{convex} if
$\alpha\prec\gamma\prec\beta$ with $\alpha,\beta\in S$ implies $\gamma\in S$. Two properties of
convex sets will be used constantly.
\begin{itemize}
\item[(C1)] The crossings of a subword $s$ of $w$ form a convex set: if
  $\alpha\prec\gamma\prec\beta$ then $\gamma$ occurs between $\alpha$ and $\beta$ in $w$.
\item[(C2)] If a convex set $S$ contains two crossings on a strand $x$, it contains all the
  crossings on $x$ between them.
\end{itemize}

\paragraph{Cyclic words.}
Theorem~\ref{thm:bands}(ii) concerns words obtained from $V_m$ by cyclic permutations as well
as commutations. To handle them, let $w$ be a braid word inducing the identity permutation,
and let $w^\infty=\cdots w\,w\,w\cdots$ be the bi-infinite periodic word. Its strands are
those of $w$, their itineraries are bi-infinite and periodic, and crossings, the relation
$\prec$ and convex sets are defined as before. The \emph{shift} $\tau$ moves each crossing
of $w^\infty$ to the corresponding crossing one period later.

\begin{lemma}\label{lem:cyclic}
Let $w$ be a braid word inducing the identity permutation, and let $w'$ be obtained from $w$
by commutations and cyclic permutations. Then the letters of $w'$ can be identified with the
crossings of a convex set $C$ of $w^\infty$, respecting strands and signs, in such a way that
$\alpha\prec\beta$ implies that $\alpha$ occurs before $\beta$ in $w'$, for
$\alpha,\beta\in C$. Consequently the crossings of any subword of $w'$ form a convex set of
$w^\infty$.
\end{lemma}

\begin{proof}
We prove more: $C$ can be chosen so that every crossing of $w^\infty$ lies in exactly one
of the sets $\tau^i(C)$, $i\in\mathbb Z$, and so that $\alpha\in\tau^i(C)$,
$\beta\in\tau^j(C)$ and $\alpha\prec\beta$ imply $i\le j$. Call the index $i$ the
\emph{level} of $\alpha$. This holds for $w'=w$, with $C$ one period of $w^\infty$; and it
implies that $C$ is convex, since a crossing between two crossings of level $0$ has
level~$0$.

If $w'$ has these properties and $w''$ is obtained from $w'$ by a commutation, the same set
$C$ serves for $w''$: the two exchanged letters $\alpha,\beta$ are adjacent in $w'$ and have no
common strand, so they are not comparable, because a chain from $\alpha$ to $\beta$ would
consist of crossings of level $0$ occurring strictly between them in $w'$.

If $w'=xy$ and $w''=yx$, let $C_x$ and $C_y$ be the sets of crossings of $x$ and of $y$,
and take $C''=C_y\cup\tau(C_x)$, with $w''$ read as $C_y$ followed by $\tau(C_x)$. Every
crossing lies in exactly one translate of $C''$. Let $\alpha\in\tau^i(C'')$ and
$\beta\in\tau^j(C'')$ with $\alpha\prec\beta$, and suppose $i>j$. The old levels of $\alpha$
and $\beta$ are $i$ or $i+1$, and $j$ or $j+1$, and the first is at most the second; so both
equal $i=j+1$, which means $\alpha\in\tau^i(C_y)$ and $\beta\in\tau^{i}(C_x)$. But
then $\beta$ occurs before $\alpha$ in the copy $\tau^i(C)$ of $w'=xy$, contradicting
$\alpha\prec\beta$. So $i\le j$. Finally, if $\alpha,\beta\in C''$ and $\alpha\prec\beta$,
then $\alpha$ occurs before $\beta$ in $w''$: this is inherited from $w'$ if both lie in
$C_y$ or both in $\tau(C_x)$, and the case $\alpha\in\tau(C_x)$, $\beta\in C_y$ is excluded by
the old levels.

The last statement follows from (C1): if $\alpha\prec\gamma\prec\beta$ with $\alpha,\beta$ in
a subword $s$ of $w'$, then $\gamma\in C$ by convexity, and $\gamma$ occurs between $\alpha$ and
$\beta$ in $w'$, hence in $s$.
\end{proof}

\subsection{The itineraries of the band braids}

Let $D_m=V_m^\infty$. We say that the bands $\mathcal P$ and $\mathcal Q$ are \emph{clasp
partners} of each other, and so are $\mathcal P'$ and $\mathcal Q'$; the bands $\mathcal P$
and $\mathcal P'$ are \emph{pass partners} of each other, and so are $\mathcal Q$ and
$\mathcal Q'$. Thus each band has one clasp partner and one pass partner, an upper band has
an upper clasp partner and a lower pass partner, and a lower band has a lower clasp partner
and an upper pass partner. The two remaining pairs of bands, $\{\mathcal P',\mathcal Q\}$
and $\{\mathcal P,\mathcal Q'\}$, are \emph{non-interacting}.

\begin{lemma}\label{lem:bandit}
In $D_m$, two strands cross only if their bands are clasp partners or pass partners. The
itinerary of each strand $x$ is a periodic sequence of \emph{blocks} of two kinds.
\begin{itemize}
\item In a \emph{clasp block}, $x$ crosses every strand of the clasp partner of its band
  exactly twice, in an order of the form $y_1,y_2,\dots,y_m,y_m,\dots,y_2,y_1$, and all these
  crossings have the same sign. Consecutive clasp blocks of $x$ have opposite signs.
\item In a \emph{pass block}, $x$ crosses every strand of the pass partner of its band exactly
  once, and all these crossings have the same sign. Pass blocks are alternately
  \emph{outward} and \emph{return} blocks; an outward block and the return block that
  follows it meet the strands in opposite orders, $z_1,\dots,z_m$ and $z_m,\dots,z_1$, and
  with opposite signs.
\end{itemize}
The blocks of an upper strand follow each other in the periodic order
\[
  \text{clasp},\ \ \text{outward},\ \ \text{return},\ \ \text{clasp},\ \ \text{outward},\ \ \text{return},
  \ \dots
\]
and those of a lower strand in the periodic order
\[
  \text{outward},\ \ \text{clasp},\ \ \text{return},\ \ \text{outward},\ \ \text{clasp},\ \ \text{return},
  \ \dots
\]
A crossing of an upper strand $x$ and a lower strand $y$ lies in an outward block of $x$ if and
only if it lies in an outward block of $y$.
\end{lemma}

\begin{proof}
This is read off from $V_m=K\cdot EFKF^{-1}E^{-1}\cdot K^{-1}\cdot EFK^{-1}F^{-1}E^{-1}$ and
from the description of $\chi_q(a,b)$ in Section~\ref{subsec:bandbraids}. In
$\chi_m(m,m)^2$, a strand of the band which starts on the left meets the strands of the other
band in the order $1,\dots,m$ during the first factor, and in the order $m,\dots,1$ during
the second, in which its band is on the right; a strand of the other band meets them in the
orders $m,\dots,1$ and then $1,\dots,m$. All crossings of $K$ are positive and all crossings of
$K^{-1}$ are negative, and $K^{-1}$ gives the same orders, since they are palindromes. In
$E$, a strand of $\mathcal P'$ meets $P_1,\dots,P_m$ in this order and a strand of
$\mathcal P$ meets $P'_m,\dots,P'_1$, with negative crossings; in $E^{-1}$ the orders are
reversed and the crossings positive. The words $F$ and $F^{-1}$ are similar, for $\mathcal Q$
and $\mathcal Q'$. The blocks of $E$ and $F$ are the outward blocks. Explicitly, the
itineraries over one period are
\[
\begin{array}{ll}
P_j: & \mathcal Q^{\,2}_{+}\ \ \mathcal P'_{-}\ \ \mathcal P'_{+}\ \ \mathcal Q^{\,2}_{-}\ \
       \mathcal P'_{-}\ \ \mathcal P'_{+}\,,\\[2pt]
Q_k: & \mathcal P^{\,2}_{+}\ \ \mathcal Q'_{+}\ \ \mathcal Q'_{-}\ \ \mathcal P^{\,2}_{-}\ \
       \mathcal Q'_{+}\ \ \mathcal Q'_{-}\,,\\[2pt]
P'_i: & \mathcal P_{-}\ \ \mathcal Q'^{\,2}_{+}\ \ \mathcal P_{+}\ \ \mathcal P_{-}\ \
        \mathcal Q'^{\,2}_{-}\ \ \mathcal P_{+}\,,\\[2pt]
Q'_l: & \mathcal Q_{+}\ \ \mathcal P'^{\,2}_{+}\ \ \mathcal Q_{-}\ \ \mathcal Q_{+}\ \
        \mathcal P'^{\,2}_{-}\ \ \mathcal Q_{-}\,,
\end{array}
\]
where $\mathcal B^{\,2}_{\pm}$ stands for a clasp block with the band $\mathcal B$,
$\mathcal B_{\pm}$ for a pass block, and the subscript is the sign.
\end{proof}

Two consequences will be used repeatedly. Let $x$ be an upper strand and $y$ a lower strand
whose bands are pass partners. Their crossings are alternately outward and return crossings.
\begin{itemize}
\item[(I1)] Between an outward crossing of $x$ and $y$ and the return crossing that follows it,
  the strand $y$ runs through a whole clasp block, and the strand $x$ crosses nothing but
  strands of the band of $y$.
\item[(I2)] Between a return crossing of $x$ and $y$ and the outward crossing that follows it,
  the strand $x$ runs through a whole clasp block, and the strand $y$ crosses nothing but
  strands of the band of $x$.
\end{itemize}

A \emph{transversal} is a set $\mathcal T$ of four strands, one from each band.

\begin{lemma}\label{lem:transversal}
For every transversal $\mathcal T$, the word $V_m|_{\mathcal T}$ is obtained from $V_1$ by
commutations, the strand of $\mathcal T$ in each band corresponding to the strand of $V_1$ in
the same band. Hence $D_m|_{\mathcal T}=D_1$. Moreover, if $S$ is a convex set of crossings of
$D_m$, the crossings of $S$ between strands of $\mathcal T$ form a convex set of
$D_m|_{\mathcal T}$.
\end{lemma}

\begin{proof}
Deleting from $\chi_q(m,m)$ all strands but one in each of the two bands leaves a single
crossing, of the same sign; so the blocks $K^{\pm1}$, $E^{\pm1}$, $F^{\pm1}$ of $V_m$ restrict
to the corresponding blocks of $V_1$. For the last statement, note that
$\alpha\prec\gamma\prec\beta$ in $D_m|_{\mathcal T}$ implies $\alpha\prec\gamma\prec\beta$ in
$D_m$, because two crossings which are consecutive on a strand of the restricted diagram
occur in the same order on this strand in $D_m$.
\end{proof}

\subsection{Pairs of strands that cross twice}

From now on, $\Sigma$ denotes a set of strands of $D_m$, and $S$ a convex set of crossings of
$D_m$ all of which are between strands of $\Sigma$. For $x,y\in\Sigma$ we write $c(x,y)$ for
the number of crossings of $x$ and $y$ in $S$. By (C2), if $S$ contains two crossings on a
strand $x$, then every strand that $x$ crosses between them belongs to $\Sigma$.

\begin{lemma}\label{lem:twice}
Suppose that $\Sigma$ contains at most $2m$ strands of the four bands. Then $c(x,y)\le2$ for all
$x,y$. If $c(x,y)=2$, then either
\begin{itemize}
\item[\textup{(a)}] the bands of $x$ and $y$ are clasp partners, and the two crossings lie in
  the same clasp block of $x$, so that they have the same sign; or
\item[\textup{(b)}] the bands of $x$ and $y$ are pass partners, the two crossings have
  opposite signs, and between them one of the strands $x,y$ runs through a whole clasp
  block, all of whose crossings belong to $S$. In particular $\Sigma$ contains a whole band,
  different from the bands of $x$ and $y$.
\end{itemize}
\end{lemma}

\begin{proof}
Let the bands of $x$ and $y$ be clasp partners, and suppose that $S$ contains two crossings of
$x$ and $y$ from different clasp blocks. By Lemma~\ref{lem:bandit}, two full pass blocks of
$x$ lie between them, and two full pass blocks of $y$. So $\Sigma$ contains the pass partners
of the bands of $x$ and of $y$. These are two different bands, containing neither $x$ nor
$y$, and $\Sigma$ would contain $2m+2$ strands. Hence all the crossings of $x$ and $y$ in $S$
lie in one clasp block, and (a) holds.

Let the bands of $x$ and $y$ be pass partners, $x$ upper and $y$ lower, and suppose that $S$
contains two crossings of $x$ and $y$; we may assume that they are consecutive among the
crossings of $x$ and $y$. If they are an outward crossing followed by a return crossing,
then by (I1) the strand $y$ runs through a clasp block between them, and $\Sigma$ contains
the clasp partner of the band of $y$. If they are a return crossing followed by an outward
crossing, then by (I2) the same holds with $x$ in place of $y$. In both cases the crossings
have opposite signs, the crossings of the clasp block belong to $S$ by (C2), and the band
obtained is not the band of $x$ or of $y$. If $S$ contained three crossings of $x$ and $y$,
both cases would occur, and $\Sigma$ would contain the clasp partners of both bands, which
are different, and $x$ and $y$: again $2m+2$ strands.
\end{proof}

\begin{lemma}\label{lem:bbcc}
Suppose that $\Sigma$ contains at most $2m$ strands of the four bands. Then there are no
strands $a,b,c\in\Sigma$ with $c(a,b)=c(a,c)=2$ such that, along $a$, both crossings with $b$
come before both crossings with $c$.
\end{lemma}

\begin{proof}
Suppose there are. By Lemma~\ref{lem:twice}, each of the pairs $\{a,b\}$ and $\{a,c\}$ is of
type~(a), a \emph{clasp pair}, or of type~(b), a \emph{handle}. Recall that all crossings on
$a$ between the first crossing with $b$ and the last crossing with $c$ belong to $S$.

\emph{Two clasp pairs.} Then $b$ and $c$ lie in the clasp partner of the band of $a$. In a
single clasp block the crossings of $a$ with $b$ and with $c$ are nested, in the order
$b,c,c,b$ or $c,b,b,c$. So the crossings with $b$ lie in one clasp block and those with $c$
in a later one. If $c$ precedes $b$ in the order of these blocks, the second crossing of $a$
with $c$ in the first block comes after the crossings with $b$, so it belongs to $S$; if $b$
precedes $c$, the first crossing of $a$ with $b$ in the later block comes before the last
crossing with $c$, so it belongs to $S$. Either way $S$ contains crossings of a clasp pair
from two clasp blocks, contradicting Lemma~\ref{lem:twice}.

\emph{Two handles.} Then $b$ and $c$ lie in the pass partner of the band of $a$; say that the
outward blocks of $a$ meet $b$ before $c$. Along $a$, the crossings with $b$ and $c$ occur in
the periodic order $b,c\,|\,c,b\,|\,b,c\,|\,c,b\,|\cdots$, the bars separating the pass
blocks. Four consecutive terms of the form $b,b,c,c$, or $c,c,b,b$, consist of a return
crossing followed by an outward crossing, for one of the two strands, and of an outward
crossing followed by a return crossing, for the other. By (I1) and (I2), $\Sigma$ contains the clasp
partner of the band of $a$ and the clasp partner of the band of $b$ and $c$. These are two
different bands, not containing $a$, and $\Sigma$ would contain $2m+1$ strands.

\emph{A clasp pair and a handle.} Say $\{a,k\}$ is the clasp pair and $\{a,z\}$ the handle,
where $\{k,z\}=\{b,c\}$. No crossing of $a$ with $k$ lies between the two crossings of $a$
with $z$. So $a$ does not run through a clasp block between these two crossings, for it
would cross $k$ there. By (I1) and (I2) it is then $z$ that runs through a clasp block between
them, and $\Sigma$ contains the clasp partner of the band of $z$; moreover, one of the two
crossings of $a$ with $z$ is an outward crossing and the other a return crossing, lying in two
consecutive pass blocks of $a$. The crossings of $a$ with $k$ lie in another block of $a$,
before or after these two pass blocks; so one of the two pass blocks lies entirely between
crossings of $S$ on $a$, and $\Sigma$ contains the whole band of $z$. Thus $\Sigma$ contains
two bands, not containing $a$, and has at least $2m+1$ strands.
\end{proof}

\subsection{Three strands}\label{subsec:three}

We need two elementary facts about braids on three strands. Recall that a word is
\emph{geodesic} if no shorter word represents the same element.

\begin{lemma}\label{lem:conjgeodesic}
Let $\{x,y\}=\{\sigma_1,\sigma_2\}$ in $B_3$, let $\varepsilon,\eta\in\{1,-1\}$, and let
$g=x^{\varepsilon}y^{2\eta}x^{-\varepsilon}$.
\begin{enumerate}
\item[\textup{(a)}] The word $g$ and all its subwords are geodesic.
\item[\textup{(b)}] $x^{\varepsilon}y^{2\eta}x^{-\varepsilon}\ne
  x^{-\varepsilon}y^{2\eta}x^{\varepsilon}$ in $B_3$.
\end{enumerate}
\end{lemma}

\begin{proof}
The automorphism of $B_3$ exchanging $\sigma_1$ and $\sigma_2$ preserves lengths, so we may
assume $x=\sigma_1$ and $y=\sigma_2$. We use the homomorphism $\rho\colon B_3\to
\mathrm{SL}_2(\mathbb Z)$ with
\[
  \rho(\sigma_1)=\begin{pmatrix}1&1\\0&1\end{pmatrix},\qquad
  \rho(\sigma_2)=\begin{pmatrix}1&0\\-1&1\end{pmatrix},\qquad\text{so that}\qquad
  \rho(g)=\begin{pmatrix}1-2\varepsilon\eta&2\eta\\-2\eta&1+2\varepsilon\eta\end{pmatrix}.
\]
Part (b) follows, since exchanging $\varepsilon$ and $-\varepsilon$ changes this matrix.
For (a), recall that two words representing the same braid have the same exponent sum, and
hence lengths of the same parity. A subword of length at most $2$ is geodesic because it is
not trivial. A subword of length $3$ is $x^{\varepsilon}y^{2\eta}$ or
$y^{2\eta}x^{-\varepsilon}$; if it were not geodesic it would equal a single letter, but both
off-diagonal entries of its image under $\rho$, which are $\pm1$ and $\pm2$, are non-zero,
unlike those of $\rho(\sigma_i^{\pm1})$. If $g$ were not geodesic it would equal a word of
length at most $2$ with exponent sum $2\eta$, that is, one of $\sigma_1^{2\eta}$,
$\sigma_2^{2\eta}$, $\sigma_1^{\eta}\sigma_2^{\eta}$, $\sigma_2^{\eta}\sigma_1^{\eta}$; each
of their images under $\rho$ has an off-diagonal entry equal to $0$ or $\pm1$, whereas those of
$\rho(g)$ are $\pm2$.
\end{proof}

\begin{lemma}\label{lem:walk}
Let $v$ be a braid word on three strands $a,b,c$ in which $b$ and $c$ do not cross, $b$ being
to the left of $c$, and in which $a$ crosses $b$ exactly $n_b\le2$ times and $c$ exactly
$n_c\le2$ times. Then the order in which $a$ meets $b$ and $c$ is determined by $n_b$, $n_c$
and the initial position of $a$, except when $a$ starts between $b$ and $c$ and $n_b=n_c=2$;
in that case the order is $b,b,c,c$ or $c,c,b,b$.
\end{lemma}

\begin{proof}
At each moment $a$ is to the left of $b$, between $b$ and $c$, or to the right of $c$; a
crossing with $b$ exchanges the first two positions and a crossing with $c$ the last two. If
$a$ starts on the left, its first crossing is with $b$, and the possible orders are
$b$; $bb$; $bc$; $bcc$; $bccb$, one for each admissible pair $(n_b,n_c)$. The case where
$a$ starts on the right is symmetric. If $a$ starts in the middle, the possible orders are
$b$; $c$; $bb$; $cc$; $bbc$; $ccb$; $bbcc$ and $ccbb$.
\end{proof}

\subsection{Proof of Theorem~\ref{thm:bands}(ii)}\label{subsec:proofii}

Let $\Omega$ be the set of words obtained from $V_m$ by commutations, cyclic permutations and
passages to the reversed inverse. A move on a word $w^{-1}$, replacing $s^{-1}$ by $s'^{-1}$,
corresponds to the move on $w$ replacing $s$ by $s'$, and one is inessential if and only if the
other is. So it is enough to consider a word $w\in\Omega$ obtained from $V_m$ by commutations
and cyclic permutations. Once we know that every move on at most $2m$ strands applicable to a
word of $\Omega$ is inessential, it follows that $\Omega$ is closed under all the operations
allowed, and all its words have length $16m^2$; this gives the last assertion of~(ii).

So let $w=p\,s\,q$, where $s$ is a $k$-strand subword with $k\le2m$, in the letters
$\sigma_a^{\pm1},\dots,\sigma_{a+k-2}^{\pm1}$, and let $s'$ be a word in the same letters with
$s'=s$ in $B_n$ and $|s'|\le|s|$. Let $\Sigma$ be the set of the $k$ strands which occupy the
positions $a,\dots,a+k-1$ at the beginning of $s$; we regard $s$ and $s'$ as braid words on
the strands of $\Sigma$. By Lemma~\ref{lem:cyclic}, the strands of $w$ are identified with
those of $D_m$, together with $n-4m$ \emph{idle} strands, which have no crossings, and the
set $S$ of crossings of $s$ is a convex set of $D_m$. So the notation and the results of the
preceding sections apply; in particular $c(x,y)$ is the number of crossings of $x$ and $y$
in $s$. We write $c'(x,y)$ for the number of crossings of $x$ and $y$ in $s'$. We must show
that $s'$ is obtained from $s$ by commutations.

\paragraph{Step 1: $c'(x,y)=c(x,y)$ for all $x,y\in\Sigma$.}
We compare $s$ and $s'$ on small sets of strands. For $T\subseteq\Sigma$ we have
$s'|_T=s|_T$ in $B_{|T|}$. Let $\mathcal T$ be a transversal. By Lemma~\ref{lem:transversal},
the crossings of $S$ between strands of $\mathcal T$ form a convex set $S_1$ of $D_1$, and
they are crossings between strands of $\mathcal T\cap\Sigma$.

\emph{Two strands.} If $\mathcal T\cap\Sigma=\{x,y\}$, then
\begin{equation}\label{eq:pair}
  c'(x,y)\ \ge\ c(x,y).
\end{equation}
Indeed, Lemma~\ref{lem:twice}, applied to $D_1$ and to the two strands $x,y$, shows that all
the crossings of $S_1$ have the same sign: case~(b) cannot occur, as it requires a third
strand. So $s|_{\{x,y\}}$ is a power $\sigma_1^{e}$ with $|e|=c(x,y)$, and the word
$s'|_{\{x,y\}}$, which represents the same element of $B_2$, has at least $|e|$ letters.

\emph{Three strands.} If $\mathcal T\cap\Sigma=\{x,y,z\}$, then
\begin{equation}\label{eq:triple}
  c'(x,y)+c'(x,z)+c'(y,z)\ \ge\ c(x,y)+c(x,z)+c(y,z).
\end{equation}
It suffices to show that the word $s|_{\{x,y,z\}}$ is geodesic in $B_3$. Among the
bands of $x,y,z$ exactly one pair is non-interacting; let $a$ be the strand of the third
band, $b$ the strand in the clasp partner of its band, and $c$ the strand in its pass partner.
In $D_1$ the strand $a$ crosses no other strands than $b$ and $c$, so by (C2) the
crossings of $S_1$ are consecutive in the itinerary of $a$ in $D_1$, which is periodic of the
form $bb,c,c,bb,c,c,\dots$ if $a$ is upper and $c,bb,c,c,bb,c,\dots$ if $a$ is lower
(Lemma~\ref{lem:bandit} for $m=1$). The set $S_1$ does not contain two crossings of $a$ and
$c$ between which $c$ runs through its clasp block, since this block consists of crossings
with the fourth strand of $\mathcal T$, which is not in $\Sigma$; by (I1) and (I2), these
are an outward crossing and the next return crossing if $a$ is upper, and a return crossing
and the next outward crossing if $a$ is lower, and in both cases they are adjacent in the
itinerary of $a$. It follows that
$S_1$ is contained in a stretch of the form $c,b,b,c$ of this itinerary, in which the two
crossings with $c$ have opposite signs and the two crossings with $b$ have equal signs. The
word of such a stretch is one of the words $g$ of Lemma~\ref{lem:conjgeodesic}, and
$s|_{\{x,y,z\}}$ is a subword of it, so it is geodesic by Lemma~\ref{lem:conjgeodesic}(a).

\emph{Counting.} At most two bands are contained in $\Sigma$, since $\Sigma$ has at most $2m$
strands. Two strands of the same band, or of non-interacting bands, do not cross in $s$, and
idle strands have no crossings in $s$. So
$|s|=\sum c(x,y)$, the sum being over the pairs of strands of $\Sigma$ from two different
bands, while $|s'|\ge\sum c'(x,y)$ over the same pairs. We show that
\begin{equation}\label{eq:count}
  \sum c'(x,y)\ \ge\ \sum c(x,y),\quad\text{with equality only if $c'(x,y)=c(x,y)$ for all
  these pairs.}
\end{equation}
Since $|s'|\le|s|$, this implies that $c'(x,y)=c(x,y)$ for these pairs and that $s'$ has no
other crossings, which is the claim of Step~1.

Suppose first that no band is contained in $\Sigma$, or that $\Sigma$ is the union of two bands.
Then for every pair $x,y\in\Sigma$ from two different bands, each of the other two bands has
a strand outside $\Sigma$, so there is a transversal with $\mathcal T\cap\Sigma=\{x,y\}$, and
\eqref{eq:count} follows from \eqref{eq:pair}.

Suppose now that $\Sigma$ contains exactly one band $\mathcal G$, and let
$Y$ be the set of the other band strands of $\Sigma$, so that $|Y|=r\le m$.
For $g\in\mathcal G$ and $y\in Y$ there is a transversal meeting $\Sigma$ in $\{g,y\}$, so
$\delta(g,y)=c'(g,y)-c(g,y)\ge0$ by \eqref{eq:pair}. For $y,z\in Y$ in different bands and
any $g\in\mathcal G$ there is a transversal meeting $\Sigma$ in $\{g,y,z\}$, and
\eqref{eq:triple} gives $c'(y,z)-c(y,z)\ge-\delta(g,y)-\delta(g,z)$. Averaging over the $m$
strands $g$ of $\mathcal G$,
\[
  c'(y,z)-c(y,z)\ \ge\ -\frac{\Delta_y+\Delta_z}{m},\qquad\text{where }\
  \Delta_y=\sum_{g\in\mathcal G}\delta(g,y)\ \ge0 .
\]
Each $y\in Y$ belongs to at most $r-1$ such pairs $\{y,z\}$. Summing over all pairs, with
$\Delta=\sum_{y\in Y}\Delta_y$,
\[
  \sum c'(x,y)-\sum c(x,y)\ \ge\ \Delta-\frac{r-1}{m}\,\Delta\ \ge\ \frac{\Delta}{m}\ \ge\ 0 .
\]
If equality holds, then $\Delta=0$, so $c'(g,y)=c(g,y)$ for all $g\in\mathcal G$, $y\in Y$;
then $c'(y,z)\ge c(y,z)$ for all pairs $y,z\in Y$, and equality of the sums forces
$c'(y,z)=c(y,z)$. This proves \eqref{eq:count}.

\paragraph{Step 2: every strand has the same itinerary in $s'$ as in $s$.}
By Step~1 two strands cross in $s'$ exactly as many times as in $s$; in particular strands
of the same band, or of non-interacting bands, do not cross in $s'$, and by
Lemma~\ref{lem:twice} two strands cross at most twice. We determine the itineraries in $s'$
in two stages.

\emph{Signs.} Let $x,y\in\Sigma$ cross in $s$. The sum of the signs of their crossings is the
same in $s'$ as in $s$, since $s'|_{\{x,y\}}=s|_{\{x,y\}}$ in $B_2$. So if the crossings of $x$ and $y$
in $s$ all have the same sign, the same holds in $s'$. Otherwise we are in case~(b) of
Lemma~\ref{lem:twice}: $x$ and $y$ cross twice, with signs $\varepsilon,-\varepsilon$ in this
order, and one of them, which we call $a$, runs between the two crossings through a whole
clasp block, contained in $S$. Let $a^*$ be the other strand and let $g$ be a strand of the
clasp partner of the band of $a$. In $s$, the strand $a$ meets $a^*$, then $g$ twice, with
equal signs $\eta$, then $a^*$ again, and it has no other crossings with $a^*$ or $g$. The bands of $a^*$ and $g$ are non-interacting, so $a^*$ and $g$ do
not cross, neither in $s$ nor in $s'$. By Lemma~\ref{lem:walk}, applied to $s|_{\{a,a^*,g\}}$
and to $s'|_{\{a,a^*,g\}}$, the strand $a$ does not start between $a^*$ and $g$, and it meets
$a^*$ and $g$ in the same order $a^*,g,g,a^*$ in $s'$. The two crossings of $a$ and $g$ have
the sign $\eta$ in $s'$ as well, by the case already treated. Hence, if the crossings of $a$ and
$a^*$ in $s'$ have the signs $\varepsilon',-\varepsilon'$ in this order, then
\[
  s|_{\{a,a^*,g\}}=\sigma_i^{\varepsilon}\sigma_j^{2\eta}\sigma_i^{-\varepsilon},\qquad
  s'|_{\{a,a^*,g\}}=\sigma_i^{\varepsilon'}\sigma_j^{2\eta}\sigma_i^{-\varepsilon'},
\]
with the same $\{i,j\}=\{1,2\}$ in both. These words represent the same braid, so
$\varepsilon'=\varepsilon$ by Lemma~\ref{lem:conjgeodesic}(b).

\emph{Orders.} Let $a\in\Sigma$ and let $b,c\in\Sigma$ be two strands which cross $a$ in $s$. The
bands of $b$ and $c$ are each a clasp partner or a pass partner of the band of $a$, so they
are equal or non-interacting; hence $b$ and $c$ do not cross, neither in $s$ nor in $s'$. By
Lemma~\ref{lem:walk}, applied to the restrictions of $s$ and $s'$ to $\{a,b,c\}$, the strand $a$ meets $b$ and $c$ in the same
order in $s'$ as in $s$, unless $a$ starts between $b$ and $c$ and crosses each of them
twice; but then $a$ would meet them in $s$ in the order $b,b,c,c$ or $c,c,b,b$, which is
excluded by Lemma~\ref{lem:bbcc}.

The itinerary of a strand $a$ is determined by the signs of its crossings with each strand
and by the relative order of its crossings with any two strands. So all strands have the same
itineraries in $s'$ as in $s$, and by Lemma~\ref{lem:itin} the word $s'$ is obtained from $s$
by commutations. This completes the proof of Theorem~\ref{thm:bands}. \qed

\begin{remark}[Where the bound $2m$ is used]
The number $2m$ enters the proof in two ways. First, a set of at most $2m$ strands cannot
contain two bands and a further strand; this is what keeps the regions of $V_m$ free of the
configurations excluded in Lemmas~\ref{lem:twice} and~\ref{lem:bbcc}. Second, in the counting
argument of Step~1, a set $\Sigma$ containing a band has at most $m$ other strands, which is
exactly what makes the coefficient $1-(r-1)/m$ positive. Both fail for $2m+1$ strands: a
strand of $\mathcal P$, the band $\mathcal P'$ and the band $\mathcal Q'$ support the essential
move used in Stage~1 of the proof of part~(iii).
\end{remark}

\begin{example}[The counting argument in a small case]\label{ex:counting}
Let $m=2$. Up to commutations, the factor $H$ of $V_2$ contains the subword
\[
  s=\sigma_2^{-1}\cdot\sigma_3\sigma_4\sigma_4\sigma_3\cdot\sigma_2 ,
\]
which involves the four strands $P'_1,P_2,Q'_1,Q'_2$, in this order at the beginning of $s$:
the strand $P'_1$ passes under $P_2$, crosses $Q'_1$ and $Q'_2$ twice each, as part of the
clasp of the lower bands, and passes back under $P_2$. Here $\Sigma$ contains the band
$\mathcal G=\mathcal Q'$, and $Y=\{P'_1,P_2\}$. There is an obvious competitor: the strand
$P_2$ lies above the other three, so by Lemma~\ref{lem:layer} it may step to the right over
$Q'_1$ and $Q'_2$ instead of stepping to the left over $P'_1$. But this costs four
crossings instead of two, and the resulting word
$\sigma_3\sigma_4\cdot\sigma_2\sigma_3\sigma_3\sigma_2\cdot\sigma_4^{-1}\sigma_3^{-1}$ has
length $8$. Step~1 shows that nothing better is possible. Let $s'$ be a word on these four
strands with $s'=s$. Taking transversals through $P'_1$ and $Q'_l$ which avoid $P_2$,
inequality \eqref{eq:pair} gives $c'(P'_1,Q'_l)\ge2$ for $l=1,2$. Taking the transversal
through $P_2$, $P'_1$ and $Q'_l$, whose fourth strand $Q_1$ or $Q_2$ is outside $\Sigma$,
inequality \eqref{eq:triple} gives
\[
  c'(P_2,P'_1)+c'(P_2,Q'_l)+c'(P'_1,Q'_l)\ \ge\ 2+0+2\qquad(l=1,2):
\]
on the three strands $P_2,P'_1,Q'_l$ the word $s$ is
$\sigma_1^{-1}\sigma_2^{2}\sigma_1$, which is geodesic. If
$c'(P'_1,Q'_l)=2$ for both $l$, adding the two inequalities gives
$2\,c'(P_2,P'_1)+c'(P_2,Q'_1)+c'(P_2,Q'_2)\ge4$, so $P_2$ takes part in at least two
crossings of $s'$, and in exactly two only if it crosses $P'_1$ twice and does not cross
$\mathcal Q'$. Hence $|s'|\ge4+2=|s|$. In general, extra crossings between $P'_1$ and $\mathcal Q'$ can be traded
against crossings of $P_2$, but only at the rate $1/m$, which is the content of the
averaging in Step~1. Step~2 then shows that, with these numbers of crossings, $s'$ is $s$.

The same subword with both strands of $\mathcal P'$, on five strands, is the one to which
the essential move of Stage~1 of the proof of part~(iii) applies: there the detour of $P_2$
over $\mathcal Q'$ costs four crossings, exactly as many as the detour over $\mathcal P'$.
\end{example}

\begin{remark}[Band braids and the family of Section~\ref{sec:family}]
By Lemma~\ref{lem:transversal}, the band braid $V_m$ contains $m^4$ copies of the hard braid
$V_1=w_{-1}(-2,2)$, one for each transversal, and the proof of Theorem~\ref{thm:bands}(ii)
reduces the rigidity of $V_m$ to simple properties of $V_1$. We expect that the same
construction, applied to the other members $w_\varepsilon(p,q)$ of the family, gives further
trivial braids with the properties of Theorem~\ref{thm:bands}; our proof uses $|p|=|q|=2$ in
Lemma~\ref{lem:walk}. Note the contrast with Section~\ref{sec:power}: for the family
$w_\varepsilon(p,q)$, what grows is the number of crossings that a three-strand move must
handle; for the band braids $V_m$, what grows is the number of strands.
\end{remark}

\section{Verification}\label{sec:verification}

That $W$ represents $1$ was checked with Artin's faithful
action of $B_4$ on the free group $F_4$ \cite{Artin1925,Artin1947}, and is consistent with the Burau
representation evaluated at random parameters. Up to cyclic permutation and
reversal, exactly $16$ words can be reached from $W$ by elementary moves, all of
length $16$, and none contains a subword $\sigma_i^{\varepsilon}\sigma_i^{-\varepsilon}$
even cyclically; so no move ever shortens the word. The list of these $16$ words is
short enough to be checked by hand. By contrast, an exhaustive search finds no
hard braids in $B_3$ of length at most $14$, in agreement with
Theorem~\ref{thm:three}.

For four and five strands we used the following exhaustive method. Every cyclically
reduced trivial word of length $2L$ can be written, after a cyclic permutation, as
$u\,v^{-1}$ with $|u|=|v|=L$ and $u=v$ in $B_n$. So grouping all freely reduced words of
length $L$ by the element they represent, computed with Artin's faithful action, and
pairing words within each group, lists every trivial cyclic word of length $2L$. The
method was checked against direct enumeration in $B_3$ at lengths $6$, $8$ and $10$.
For each word obtained, the set of words reachable from it by elementary moves was
searched, up to a cap of $10^5$ words. In $B_4$ there are $39\,536$ cyclically reduced
trivial words of length $14$, up to cyclic permutation and reversed inverse; the cap was
never reached, and none of these words, nor any of length $10$ or $12$, is hard. Since a
trivial word has exponent sum $0$ and hence even length, $W$ is a shortest hard braid in
$B_4$. Among the $483\,625$ such words of length $16$, the search found exactly $16$
hard ones, all of them members of the family of Section~\ref{sec:family}; we did not
record whether the cap was reached at this length. In $B_5$ the same method finds no hard
braids of length at most $12$, and the cap was never reached.

For the family of Section~\ref{sec:family}, every step of the proof of
Lemma~\ref{lem:slide} was checked by machine to be a legal
length-preserving move, for $j=2,3,4$ and $-7\le q\le7$. Independently, the
exhaustive search of all words reachable by elementary moves confirmed that
$w_\varepsilon(p,q)$ is hard in all $490$ cases with opposite signs, $p,q\in\{\pm2,\pm3,\pm4,5\}$
and $(n,j)\in\{(4,2),(5,2),(5,3),(6,3),(6,4)\}$, in agreement with
Theorem~\ref{thm:hard}. Triviality was checked with Artin's faithful action of $B_n$
on the free group.

For the braids $W_k$ of Proposition~\ref{prop:exponent} with $1\le k\le4$, a search of
all words reachable by elementary moves and power braid moves of exponent at most
$2k-1$ finds the same $16$ words as for $W$, none admitting a deletion, while allowing
exponent $2k$ makes a deletion reachable. For $1\le k\le3$, a complete untangling of
$W_k$ with power braid moves of exponent at most $2k$, never increasing the length, was
also found by computer.

\paragraph{Moves on three strands.}
Before the band braids of Section~\ref{sec:bands} were found, we tested statement
$(\mathrm T_3)$, that two-generator moves always suffice, in two ways.

First, on every hard braid we had found. Besides the exhaustive searches above, we searched
trivial commutators $[\sigma_j^{-p},y]$ in $B_4$ and $B_5$, with $y$ ranging over short
conjugates of powers of generators and over all words of length at most $7$ in $B_4$
commuting with $\sigma_j$, as well as trivial words $u\,v^{-1}$ in which $v$ is a random
respelling of a random word $u$. The latter were never hard. Up to cyclic permutation and
reversed inverse, these searches produced $50$ hard braid words, of lengths $16$ to $22$, and
all of them are carried to the empty word, without increasing their length, both by
two-generator moves and by elementary and power braid moves. Two-generator moves on subwords
in $\sigma_i,\sigma_{i+1}$ were implemented with a table of all freely reduced words of length
at most $10$ in $B_3$, sorted by the element they represent. All $50$ words, however, belong
to the family of Section~\ref{sec:family}, to which Theorem~\ref{thm:power} already applies;
so these tests did not put the statement under strain. For $1\le k\le4$, two-generator moves
on subwords of length at most $2k+1$ never make a shortening of $W_k$ reachable, while
subwords of length $2k+2$ do, in agreement with Proposition~\ref{prop:boundedlength}.

Second, exhaustively on short words. For every trivial word of a given length we computed
the set of words reachable from it by commutations and by length-preserving moves on three
strands, that is, by replacing a maximal subword in two adjacent generators by any word of
the same length representing the same element of $B_3$; and we checked whether some word of
this set contains a subword in two adjacent generators which is not geodesic in $B_3$, and
can therefore be shortened. This is the case for every freely reduced trivial word of length
at most $16$ in $B_4$ (about $2.3\cdot10^{7}$ words of length $16$), at most $14$ in $B_5$
(about $7.2\cdot10^{7}$ words of length $14$) and at most $12$ in $B_6$ (about
$3.1\cdot10^{7}$ words of length $12$). By induction on the length, every trivial word within
these bounds can be untangled monotonically by moves on at most three strands, even without
cyclic permutations. So a counterexample to $(\mathrm T_3)$ cannot be short; the one we know,
$V_2$, has $64$ crossings on $8$ strands.

\paragraph{Band braids.}
Theorem~\ref{thm:bands} is proved by hand, but we also checked it by machine in small cases.
Triviality of $V_m$ was checked with Artin's faithful action for $m\le3$, and so was the fact
that the restriction of $V_m$ to each of the $m^4$ transversals is $V_1$ up to commutations
(Lemma~\ref{lem:transversal}), for $m=2,3$.

For part~(ii) we used the following method, which is independent of the proof in
Section~\ref{sec:bands}. In the periodic word $V_m^\infty$, the sets of crossings which can be
made into a $k$-strand subword by commutations and cyclic permutations are the convex sets
of crossings whose letters lie in $k-1$ consecutive generators; the maximal ones are the
maximal cliques of a compatibility graph, in which two such crossings are joined unless a
crossing with a letter outside these generators lies between them. If a move applies to a
subword, a move applies to the maximal convex set containing it. So it is enough to list the
maximal convex sets and to check that each of them is \emph{rigid}: every word in the same
generators, representing the same braid and no longer than it, is obtained from it by
commutations. Rigidity was tested by listing all such words, with a meet-in-the-middle
search and Artin's faithful action.
\begin{itemize}
\item For $m=2$ and $k=4$ there are $33$ distinct maximal words, of length at most $8$, and all
  of them are rigid in $B_4$. This proves Theorem~\ref{thm:bands}(ii) for $V_2$ by machine;
  in particular it confirms, independently of Section~\ref{sec:bands}, that $V_2$ is a
  counterexample to $(\mathrm T_3)$. Two further programs, one computing the maximal convex
  sets by a different method and one testing all maximal subwords in three consecutive
  generators of $20\,000$ words obtained from $V_2$ by random commutations, gave the same
  result.
\item For $m=3$ and $k=6$ there are $120$ distinct maximal words, of length at most $18$.
  The $110$ words of length at most $10$ are rigid in $B_6$, by the same exhaustive search.

  Eight words have length $16$; in them a strand of one band steps over two strands of
  another band while these form a clasp with a third band. For words of this length the
  search with Artin's action is too large, and we used a fingerprint instead: the image of a
  word under the Burau representation, with the variable specialised to a random element of
  a large prime field. Equal braids have equal fingerprints, so the number of words with the
  fingerprint of $s$ is an upper bound for the number of words representing the same braid
  as $s$; that the Burau representation is not faithful on six strands does not matter for
  this purpose. For each of the eight words $s$, a meet-in-the-middle count over all pairs
  of words of length $8$ finds exactly $196$ words of length $16$ with the fingerprint of
  $s$, and no word of length $14,12,\dots,0$. Since $196$ is also the number of words
  obtained from $s$ by commutations, all of which represent the same braid as $s$, there
  are no others, and $s$ is rigid.

  The last two words, of length $18$, are the clasps $K^{\pm1}$. A word no longer than $K$
  and representing the same braid has the same exponent sum $18$, so it is a positive word;
  by Garside's theorem \cite{Garside} it is obtained from $K$ by positive braid relations
  and commutations; and we checked that no braid relation applies to any word obtained from
  $K$ by commutations. So $K^{\pm1}$ are rigid, and Theorem~\ref{thm:bands}(ii) is confirmed
  by machine for $V_3$ as well.
\item For $m=2,3$ we also checked the statements of Lemmas~\ref{lem:twice}
  and~\ref{lem:bbcc} on all maximal convex sets.
\end{itemize}
For part~(iii), the untangling of $V_m$ described in Section~\ref{subsec:layers} was
carried out by machine for $m=1,2,3$: each of the $4m$ moves was checked to replace a subword
on $2m+1$ adjacent strands by an equal braid of the same length, and commutations and
deletions then lead to the empty word.

\section{Conclusion and open questions}\label{sec:conclusion}

We have shown four things. First, hard braids exist in every braid group $B_n$ with
$n\ge4$, and there are many of them: the family $w_\varepsilon(p,q)$ of
Theorem~\ref{thm:hard} has members of every even length from $16$ on. Second, this is
optimal in the number of strands: by Theorem~\ref{thm:three} there are no hard braids on
three strands, where the elementary moves always suffice to untangle a trivial braid
without increasing its length; and an exhaustive computer search shows that our shortest
example $W$ is a shortest hard braid on four strands. Third, the hard braids of this family
are not very hard: by Theorem~\ref{thm:power}, power braid moves untangle all of them
monotonically, although the exponent needed is unbounded, by
Proposition~\ref{prop:exponent}. Fourth, this is not the general picture: by
Theorem~\ref{thm:bands}, the band braid $V_m$ on $4m$ strands is not affected by any move on
at most $2m$ strands, and is untangled by moves on $2m+1$ strands.

Four nested statements were in view:
\begin{enumerate}
\item every trivial braid word can be untangled monotonically by elementary moves;
\item every trivial braid word can be untangled monotonically by elementary moves and
  power braid moves;
\item every trivial braid word can be untangled monotonically by two-generator moves, that
  is, by moves on at most three strands;
\item for some $k$, every trivial braid word can be untangled monotonically by moves on at
  most $k$ strands.
\end{enumerate}
Each statement implies the next, and all four are false (Table~\ref{tab:summary}): the first by
Theorem~\ref{thm:hard}, and the fourth, hence all of them, by Theorem~\ref{thm:bands}. The
family of Section~\ref{sec:family} separates the first statement from the second, in the
sense that it refutes the first and satisfies the second; the band braid $V_2$ refutes the
second and the third and satisfies the fourth with $k=5$. In the terms of the introduction:
on three strands braids behave like words in Coxeter groups, even for the elementary moves;
on many strands they behave like diagrams of knots, for every notion of local move that
bounds the number of strands involved.

\begin{table}[ht]
\centering
\renewcommand{\arraystretch}{1.25}
\begin{tabular}{p{0.30\textwidth}p{0.30\textwidth}p{0.30\textwidth}}
\hline
Moves allowed & They untangle & They do not untangle\\
\hline
elementary moves & every trivial word on $3$ strands (Theorem~\ref{thm:three}) &
  $w_\varepsilon(p,q)$, on $4$ strands (Theorem~\ref{thm:hard})\\
elementary moves and power braid moves of exponent at most $M$ &
  $w_\varepsilon(p,q)$ with $|q|\le M$ (Theorem~\ref{thm:power}) &
  $w_\varepsilon(p,q)$ with $|p|,|q|>M$ (Proposition~\ref{prop:exponent})\\
two-generator moves, that is, moves on at most $3$ strands &
  every $w_\varepsilon(p,q)$; every trivial word of length at most $16$ on $4$ strands &
  $V_2$, on $8$ strands (Theorem~\ref{thm:bands})\\
moves on at most $2m$ strands & & $V_m$, on $4m$ strands (Theorem~\ref{thm:bands})\\
moves on at most $2m+1$ strands & $V_m$ (Theorem~\ref{thm:bands}) & $V_{m+1}$\\
\hline
\end{tabular}
\caption{Which moves untangle which trivial braid words monotonically.}
\label{tab:summary}
\end{table}

These results leave the following questions open.
\begin{enumerate}
\item What is the length of a shortest hard braid in $B_n$ for $n\ge5$? In $B_4$ it is
  $16$, attained by $W$. In $B_5$ there are no hard braids of length at most $12$, so the
  answer is $14$ or $16$ for $n=5$.
\item What is the least number of strands, and the least length, of a trivial braid word
  that cannot be untangled monotonically by two-generator moves? Theorem~\ref{thm:bands}
  gives an example with $64$ crossings on $8$ strands, and the search of
  Section~\ref{sec:verification} shows that there is none of length at most $16$, $14$,
  $12$ on $4$, $5$, $6$ strands respectively. In particular, is statement $(\mathrm T_3)$
  true in $B_4$? More generally, what is the least number $\kappa(n)$ of strands such that
  moves on at most $\kappa(n)$ strands untangle every trivial braid word on $n$ strands?
  Corollary~\ref{cor:tits} gives $2\lfloor n/4\rfloor+1\le\kappa(n)\le n$; clearly
  $\kappa(2)=2$ and $\kappa(3)=3$, and $\kappa(4)$ is $3$ or $4$.
\item For the elementary moves alone, how many extra crossings must be inserted to
  untangle a hard braid? Does this number grow along the family $W_k$? The same question
  can be asked for the band braids $V_m$ and moves on at most $2m$ strands.
\item Is there a finite set of length-non-increasing moves that untangles all the braids
  $W_k$ of Section~\ref{sec:power}? By Propositions~\ref{prop:exponent}
  and~\ref{prop:boundedlength}, power braid moves of bounded exponent, and two-generator
  moves of bounded length, do not.
\item Are there hard braids of yet another nature? In $B_4$, every hard braid of length $16$
  found by our exhaustive search belongs to the family of Section~\ref{sec:family}, and the
  band braids are obtained from a member of this family by replacing strands by bands. Do
  the bands obtained from the other members $w_\varepsilon(p,q)$ have the properties of
  Theorem~\ref{thm:bands}?
\item Dynnikov's theorem \cite{Dynnikov2006} shows that the unknot can be simplified
  monotonically if one uses arc-presentations and exchange moves, which are not confined to
  a bounded part of the diagram. Is there a natural family of moves on braid words, necessarily
  involving unboundedly many strands, and a natural notion of complexity, for which every
  trivial braid can be simplified monotonically?
\end{enumerate}

\section{AI Declaration}

The author has published scores of academic papers, and for many of them, he wrote computer code searching for examples or counterexamples. This is the first article that the author would not be able to complete without AI. The first example of a hard braid was found after several weeks of giving prompts to Claude AI software (version Opus 5.5), with the exchange of ideas proceeding along the lines of what is, in this paper, Section \ref{sec:diagram} and Figure \ref{fig:basketW}. Then to produce a more complicated example, the braid in Figure \ref{fig:band} in Section \ref{sec:bands}, the author used Claude AI software of a more advanced version Fable 5.2.

\end{document}